\documentclass[10pt,reqno]{amsart}
\usepackage{amssymb}
\usepackage{amsmath}
\input xy
\xyoption{all}

\newtheorem{theorem}{Theorem}[section]

\theoremstyle{definition}

\theoremstyle{remark}
\newtheorem{remark}[theorem]{Remark}

\numberwithin{equation}{section}

\newcommand{\nc}{\newcommand}
\nc{\cal}{\mathcal} %% Needed for LaTeX2e
\nc{\la}{\langle} \nc{\ra}{\rangle}
 \nc{\CA}{\cal A}
 \nc{\CBB}{\cal B}
 \nc{\CC}{\cal C}
\nc{\CDD}{\cal D} \nc{\CE}{\cal E} \nc{\CF}{\cal F} \nc{\CG}{\cal
G} \nc{\CH}{\cal H} \nc{\CI}{\cal I} \nc{\CJ}{\cal J}
\nc{\CK}{\cal K} \nc{\CL}{\cal L} \nc{\CM}{\cal M} \nc{\CN}{\cal
N} \nc{\CO}{\cal O} \nc{\CP}{\cal P} \nc{\CQ}{\cal Q}
\nc{\CR}{\cal R} \nc{\CS}{\cal S} \nc{\CT}{\cal T} \nc{\CU}{\cal
U} \nc{\CV}{\cal V} \nc{\CW}{\cal W} \nc{\CZ}{\cal Z}

\nc{\fa}{\mathfrak a} \nc{\fg}{\mathfrak g} \nc{\fk}{\mathfrak k}
\nc{\fh}{\mathfrak h} \nc{\fm}{\mathfrak m} \nc{\fn}{\mathfrak n}
\nc{\fA}{\mathfrak A} \nc{\fC}{\mathfrak C} \nc{\fI}{\mathfrak I}
\nc{\fL}{\mathfrak L} \nc{\fS}{\mathfrak S}

\nc{\nen}{\newenvironment} \nc{\ol}{\overline}
\nc{\ul}{\underline} \nc{\lra}{\longrightarrow}
\nc{\lla}{\longleftarrow} \nc{\Lra}{\Longrightarrow}
\nc{\Lla}{\Longleftarrow} \nc{\Llra}{\Longleftrightarrow}
\nc{\hra}{\hookrightarrow} \nc{\iso}{\overset{\sim}{\lra}}
\nc{\Hom}{\mathrm{Hom}} \nc{\Mor}{\mathrm{Mor}}

\nc{\notebox}[1]{\noindent\fbox{\parbox{12.5cm}{\sf #1}}\\[8pt]}
\nc{\Thm}[1]{Theorem~\ref{#1}} \nc{\Prop}[1]{Proposition~\ref{#1}}
\nc{\Lem}[1]{Lemma~\ref{#1}} \nc{\Cor}[1]{Corollary~\ref{#1}}
\nc{\Conj}[1]{Conjecture~\ref{#1}} \nc{\Claim}[1]{Claim~\ref{#1}}
\nc{\Defn}[1]{ Definition~\ref{#1}} \nc{\Exa}[1]{Example~\ref{#1}}
\nc{\Rem}[1]{Remark~\ref{#1}} \nc{\Note}[1]{Note~\ref{#1}}

\marginparsep 0.1cm \marginparwidth 2.5cm \nc{\marg}{\marginpar}

\nen{thm}[1]{\label{#1}{\bf Theorem.} \em}{}
\nen{thma}[1]{\label{#1}{\bf Theorem A.\ } \em}{}
\nen{thmb}[1]{\label{#1}{\bf Theorem B.\ } \em}{}
\nen{prop}[1]{\label{#1}{\bf Proposition.\ } \em}{}
\nen{lem}[1]{\label{#1}{\bf Lemma.\ } \em}{}
\nen{klem}[1]{\label{#1}{\bf Key Lemma.\ } \em}{}
\nen{cor}[1]{\label{#1}{\bf Corollary.\ } \em}{}
\nen{cora}[1]{\label{#1}{\bf Corollary A.\ } \em}{}
\nen{corb}[1]{\label{#1}{\bf Corollary B.\ } \em}{}
\nen{conj}[1]{\label{#1}{\bf Conjecture.\ } \em}{}
\nen{claim}[1]{\label{#1}{\bf Claim.\ } \em}{}

\nen{defn}[1]{\label{#1}{\bf Definition.\ } }{}
\nen{exa}[1]{\label{#1}{\bf Example.\ } }{}

\nen{rem}[1]{\label{#1}{\em Remark.\ } }{}
\nen{exer}[1]{\label{#1}{\em Exercise.\ } }{}
\nen{pf}{\begin{proof}}{\end{proof}}

 \nc{\br}{\mathbb R}
 \nc{\bz}{\mathbb Z}
 \nc{\bc}{\mathbb C}
 \nc{\bn}{\mathbb N}
 \nc{\ck}{\mathcal{K}}
 \nc{\G}{\Gamma}
 \nc{\sm}{\setminus}
 \nc{\sub}{\subset}
 \nc{\lm}{\lambda}
 \nc{\al}{\alpha}
 \nc{\bt}{\beta}
 \nc{\om}{\omega}
 \nc{\dl}{\delta}
 \nc{\g}{\gamma}
 \nc{\Dl}{\Delta}
 \nc{\Om}{\Omega}
 \nc{\s}{\sigma}
 \nc{\ro}{\rho}
 \nc{\te}{\theta}
 \nc{\SLR}{SL_2(\br)}
 \nc{\GLR}{GL_2(\br)}
 \nc{\PGLR}{PGL_2(\br)}
 \nc{\PSLR}{PSL_2(\br)}
 \nc{\SLC}{SL(2,\bc)}
 \nc{\uH}{\mathbb H}
 \nc{\fD}{\mathcal{D}}
 \nc{\fE}{\mathcal{E}}
 \nc{\fO}{\mathcal{O}}
 \nc{\haf}{\frac{1}{2}}
 \nc{\qtr}{\frac{1}{4}}
 \nc{\shaf}{{\scriptstyle\frac{1}{2}}}
 \nc{\hlm}{{\scriptstyle\frac{\lambda}{2}}}
 \nc{\E}{\mathbb{E}}

 \nc{\8}{\infty}
 \nc{\7}{{-\infty}}
 \nc{\inv}{^{-1}}
 \nc{\eps}{\varepsilon}

\begin{document}

\title[Rankin-Selberg without unfolding]
{Rankin-Selberg without unfolding and bounds for spherical Fourier
coefficients of Maass forms}

\author{Andre Reznikov}
\thanks{The research was partially supported  by a BSF grant, by the
Minerva Foundation, by the Excellency Center ``Group Theoretic
Methods in the Study of Algebraic Varieties'' of the  Israel
Science Foundation and the Emmy Noether Institute for Mathematics
(the Center of Minerva Foundation of Germany). The paper was
mostly written during one of the visits of the author to MPIM at
Bonn. It is a pleasure to thank MPIM for its excellent working
atmosphere.}
\address{Bar-Ilan University, Ramat-Gan, Israel}
\email{reznikov@math.biu.ac.il}

% General info

%\date{\today}

\dedicatory{To Joseph Bernstein, as a small token of gratitude.}


\begin{abstract} We use the uniqueness of various invariant
functionals on irreducible unitary representations of $\PGLR$ in
order to deduce the classical Rankin-Selberg identity for the sum
of Fourier coefficients of Maass cusp forms and its new
anisotropic analog. We deduce from these formulas non-trivial
bounds for the corresponding unipotent and spherical Fourier
coefficients of Maass forms. As an application we obtain a
subconvexity bound for certain $L$-functions. Our main tool is the
notion of a Gelfand pair from the representation theory.
\end{abstract}

\maketitle
\section{Introduction}
\label{intro}

In this paper we study periods of automorphic functions. We
present a new method which allows one to obtain non-trivial
spectral identities for weighted sums of certain periods of
automorphic functions. These identities are modelled on the
classical identity of R. Rankin \cite{Ra} and A. Selberg
\cite{Se}. We recall that the Rankin-Selberg identity relates the
weighted sum of Fourier coefficients of a cusp form $\phi$ to the
weighted integral of the inner product of $\phi^2$ with the
Eisenstein series (e.g., formula \eqref{RS-Fourier} below).

We show how to deduce the classical Rankin-Selberg identity and
similar new identities from the uniqueness principle in
representation theory (also known under the following names: the
multiplicity one property, Gelfand pair). The uniqueness principle
is a powerful tool in representation theory; it plays an important
role in the theory of automorphic functions.

We associate a non-trivial spectral identity  to certain {\it
pairs} of different {\it triples} of Gelfand subgroups. Namely, we
associate a spectral identity (see the formula
\eqref{spec-two-trip} below) with two triples
$\mathcal{F}\subset\mathcal{H}_1\subset\mathcal{G}$ and
$\mathcal{F}\subset \mathcal{H}_2\subset\mathcal{G}$ of subgroups
in a group $\mathcal{G}$ such that pairs
$(\mathcal{G},\mathcal{H}_i)$ and $(\mathcal{H}_i,\mathcal{F})$
for $i=1,\ 2,$ are strong Gelfand pairs having the same subgroup
$\mathcal{F}$ in the intersection (for the notion of Gelfand pair
that we use, see Section \ref{Gelf-period-sect}). We call such a
collection $(\mathcal{G},\mathcal{H}_1,
\mathcal{H}_2,\mathcal{F})$ a {\it strong Gelfand formation}.

In the reminder of the Introduction we explain our general idea
and describe how to implement it in order to reprove the classical
Rankin-Selberg formula. We also obtain a new anisotropic analog of
the Rankin-Selberg formula. We present then an analytical
applications of these spectral identities towards non-trivial
bounds for various Fourier coefficients of cusp forms. The novelty
of our results lies mainly  in the method, as we do not rely on
the well-known technique of Rankin and Selberg which is called
``unfolding''. Instead, we use the uniqueness of relevant
invariant functionals which we explain below.

We would like to mention that the Rankin-Selberg method was
revolutionized by H. Jacquet, J. Shalika and I. Piatetski-Shapiro
who constructed integral representations for many automorphic
$L$-functions (see the survey paper \cite{Bu} and references
therein). It is not yet clear what is the relation between these
works and our point of view.

\subsection{Periods and Gelfand pairs} We
briefly review some well-known notions and constructions from
representation theory and the theory of automorphic functions
needed to formulate our identities.

Taking periods is the classical technique to study automorphic
functions. It goes back, at least, to Hecke and Maass. In the
classical language it means the following. Let $\phi$ be an
automorphic function on $Y$ and let $\psi$ be an automorphic
function on a cycle $C\subset Y$, where the cycle $C$ is equipped
with a measure $dc$. Then one can consider the period defined via
the integral
\begin{equation*}
\int_C \phi(c)\psi(c)\ dc\ .
\end{equation*}

It is well-known that in the modern language of automorphic
representations this construction leads to the following setup
which we are going to use throughout the paper.

\subsubsection{Automorphic representations}
Let $\mathcal{G}$ be a real  reductive Lie group (e.g., $\SLR$),
let $\G_\mathcal{G}\subset \mathcal{G}$ be a lattice and let
$X_\mathcal{G}=\G_\mathcal{G}\sm \mathcal{G}$ be the corresponding
automorphic space equipped with a $\mathcal{G}$-invariant measure
$\mu_\mathcal{G}$ (which we will always normalize to have the
total mass one). We denote by
$L^2(X_\mathcal{G})=L^2(X_\mathcal{G},\mu_\mathcal{G})$ the
corresponding unitary representation of $\mathcal{G}$ and by
$C^\8(X_\mathcal{G})$ its smooth part.

Let $(\pi, \mathcal{G}, V)$ be an abstract unitary irreducible
representation of $\mathcal{G}$. We will work solely with the
spaces of smooth vectors. Hence by $(\pi, V)$ we mean the
representation in the space of smooth vectors equipped with the
invariant Hermitian form. (Usually, one denotes a unitary
representation by $(\pi, \mathcal{G}, L)$ where $L$ is the
corresponding Hilbert space and the space of smooth vectors is
denoted by $V=L^\8$.)

A pair $(\pi,\nu)$, where $\nu:V\to L^2(X_\mathcal{G})$ is an
isometric $\mathcal{G}$-morphism, is called an {\it automorphic}
representation of $\mathcal{G}$ on $X_\mathcal{G}$. It is known
that $\nu: V\to C^\8(X_\mathcal{G})$ and we will denote by
$V^{aut}_\pi=V_{(\pi,\nu)}\subset C^\8(X_\mathcal{G})$ the image
of $V$ under $\nu$. We denote by $\nu^*:C^\8(X_\mathcal{G})\to V$
the adjoint map.

\subsubsection{Periods and representation theory} The notion of a
period has the following counterpart in the language of
automorphic representations.

Let $\mathcal{H}\subset \mathcal{G}$ be a subgroup (not
necessarily reductive, e.g., a unipotent subgroup of $\SLR$) and
let $X_\mathcal{H}\subset X_\mathcal{G}$ be a closed orbit of
$\mathcal{H}$ (e.g.,
$X_\mathcal{H}=\mathcal{H}\cap\G_\mathcal{G}\sm \mathcal{H}$). We
fix a $\mathcal{H}$-invariant measure $\mu_\mathcal{H}$ on
$X_\mathcal{H}$ and consider the unitary representation
$L^2(X_\mathcal{H},\mu_\mathcal{H})$ of $\mathcal{H}$. We denote
by $r_\mathcal{H}=r_{X_\mathcal{H}}: C^\8(X_\mathcal{G})\to
C^\8(X_\mathcal{H})$ the corresponding restriction map.

Let $(\pi,\nu_\pi)$ be an automorphic representation of
$\mathcal{G}$ on $X_\mathcal{G}$ and let $(\s,\nu_\s)$ an
automorphic representation of $\mathcal{H}$ on $X_\mathcal{H}$.
This means that for the abstract irreducible unitary
representation $\s$ in the space of smooth vectors $W$ equipped
with the $\mathcal{H}$-invariant Hermitian form, we fix an
isometry $\nu_\s:W\to L^2(X_\mathcal{H})$; we allow ourselves to
use the name automorphic even if $\mathcal{H}$ is not reductive.
Let $\nu^*:C^\8(X_\mathcal{H})\to W$ be the adjoint map (induced
by the $\mathcal{H}$-invariant Hermitian form on
$L^2(X_\mathcal{H},\mu_\mathcal{H})$).

We consider the $\mathcal{H}$-{\it equivariant} map
$T^{aut}_{\pi,\s}=T^{aut}_{X_\mathcal{H},\pi,\s}=\nu^*_\s\circ
r_{\mathcal{H}}\circ\nu_\pi:V\to W$ defined via the composition
\begin{equation*}
 T^{aut}_{\pi,\s}:V\stackrel{\nu_\pi}{\longrightarrow}C^\8(X_\mathcal{G})
 \stackrel{r_{\mathcal{H}}}{\longrightarrow}C^\8(X_\mathcal{H})
 \stackrel{\nu_{\s}^*}{\longrightarrow}W\ .
\end{equation*}
We call this map the automorphic period map (or simply the period)
associated to the collection $(\pi,\nu_\pi)$, $(\s,\nu_\s)$,
$X_\mathcal{G}$, $X_\mathcal{H}$ and the choice of corresponding
measures. This is the representation theoretic substitute for the
classical period. It is well-defined if $X_\mathcal{H}$ is
compact. Otherwise we have to assume that $\pi$ or $\s$ is
cuspidal (i.e., rapidly decaying at infinity if $X_\mathcal{H}$ is
not compact). For a non-cuspidal data one has to introduce
appropriate regularization in order to make sense out of the
automorphic period.

Clearly $T^{aut}_{\pi,\s}\in\Hom_\mathcal{H}(V,W)$. We denote the
vector space $\Hom_\mathcal{H}(V,W)$ by $\mathcal{P}(V,W)$ and
call it the vector space of {\it periods} between $\pi$ and $\s$.
The space of periods is defined even if $\pi$ or $\s$ are not
automorphic.

\subsubsection{Gelfand pairs and periods}\label{Gelf-period-sect}
In many cases interesting periods are associated to some special
subgroup (or representations). Many of these examples come from
what is called the multiplicity one representations or Gelfand
pairs (see \cite{Gr} and references therein). In what follows we
will use the notion of Gelfand pairs for {\it real Lie groups}.

A pair $(\mathcal{G},\mathcal{H})$ of a real Lie group
$\mathcal{G}$ and a real Lie subgroup $\mathcal{H}\subset
\mathcal{G}$ is called a {\it strong Gelfand pair} if for any pair
of irreducible representations $(\pi,V)$ of $\mathcal{G}$ and
$(\s,W)$ of $\mathcal{H}$, the dimension of the space
$\Hom_\mathcal{H}(V, W)$ of $\mathcal{H}$-equivariant maps (i.e.,
the space of periods $\mathcal{P}(V,W)$) for the spaces of {\it
smooth} vectors is at most {\it one}. (It is well-known in
representation theory of reductive groups that the dimension of
the space $\Hom_\mathcal{H}(V,W)$ for {\it infinite-dimensional}
representations plays the role of the index for {\it
finite-dimensional} representations.) The pair is called a Gelfand
pair if the same holds for $\s$ being the trivial representation.

One of the important observations in the theory of automorphic
functions is that for Gelfand pairs, automorphic periods
$T^{aut}_{\pi,\s}$ lead to certain interesting numbers or
functions (e.g., Fourier coefficients, $L$-functions). This is
based on the fact that in many cases one can construct another
vector in the {\it one-dimensional} vector space
$\mathcal{P}(V,W)$.

Namely, usually representations $\pi$ and $\s$ could be
constructed explicitly in some {\it model} spaces of sections of
various vector bundles over appropriate homogeneous manifolds
(e.g., for $\SLR$, in the spaces of homogeneous functions on the
plane $\br^2\sm 0$; see Section \ref{reps-sect}). These models
usually exist for all representations and not only for the
automorphic ones. If $\dim\mathcal{P}(V,W)=1$, using these models
one can construct an explicit non-zero map $T^{mod}_{\pi,\s}
\in\mathcal{P}(V,W)$ by means of the corresponding kernel (i.e.,
construct this map as an integral operator with an explicit
kernel). We call such a map a {\it model} period.

When $\dim\mathcal{P}(V,W)=1$, such a choice of a non-zero model
period $T^{mod}_{\pi,\s}$ gives rise to the automorphic
coefficient of proportionality
$a_{\pi,\s}=a_{\nu_\pi,\nu_\s}\in\bc$ such that
\begin{equation*}
T^{aut}_{\pi,\s}=a_{\pi,\s}\cdot T^{mod}_{\pi,\s}\ .
\end{equation*}
We would like to study these constants. In many cases these
constants are related to interesting objects (e.g., Fourier
coefficients of cusp forms, special values of $L$-functions etc.).
Of course, these constants depend on the choice of model periods
$T^{mod}_{\pi,\s}$, normalization of measures, etc. In many cases
we can find a way to canonically normalize norms of model maps in
the ad\`elic setting (and hence define canonically if not the
constants themselves then their absolute values). We hope to
discuss this question elsewhere.

We explain now how sometimes one can obtain spectral identities
involving certain coefficients $a_{\pi,\s}$.

\subsubsection{Triples of Gelfand subgroups and
periods}\label{tripl-Gelf-sect} Let $\mathcal{G}$ be a real Lie
group and $\mathcal{F}\subset\mathcal{G}$ a real Lie subgroup
which is {\it not} a Gelfand subgroup. Suppose we still want to
study a period of an automorphic representation $\pi$ of
$\mathcal{G}$ with respect to an $\mathcal{F}$-invariant closed
cycle $X_\mathcal{F}\subset X_\mathcal{G}$ (endowed with a
$\mathcal{F}$-invariant measure $\mu_\mathcal{F}$) and some
automorphic representation $\chi$ of $\mathcal{F}$. For example,
for the trivial representation of $\mathcal{F}$ we obtain  the
$\mathcal{F}$-invariant functional $I_\mathcal{F}$ on $\pi$ given
by the integral over the cycle $X_\mathcal{F}$:
\begin{equation*}\label{per-F}
I_\mathcal{F}(v)=\int_{X_\mathcal{F}}\nu_\pi(v)d\mu_\mathcal{F}\ \
{\rm for\ any\ }v\in V.
\end{equation*}
We can not apply the idea described above directly. Instead in
some cases we can obtain a spectral decomposition of this period.

Namely, suppose we can find an intermediate subgroup
$\mathcal{H}$, $\mathcal{F}\subset\mathcal{H}\subset\mathcal{G}$,
and an intermediate closed cycle $X_\mathcal{F}\subset
X_\mathcal{H}\subset X_\mathcal{G}$ such that both pairs
$(\mathcal{G},\mathcal{H})$ and $(\mathcal{H},\mathcal{F})$ are
strong Gelfand pairs. We claim that this leads to the spectral
decomposition of the functional $I_\mathcal{F}$. In what follows
we discuss for simplicity only the case of the trivial
representation of $\mathcal{F}$. The case of a non-trivial
representation $\chi$ is similar (and leads to other interesting
identities).

Consider the space
$L^2(X_\mathcal{H})=L^2(X_\mathcal{H},\mu_\mathcal{H})$ and let
assume that it has a decomposition into a direct sum (in general a
direct integral)
\begin{equation}\label{spec-decomp-H}
L^2(X_\mathcal{H})=\oplus_i\s_i
\end{equation}
of irreducible automorphic representations
$(\s_i,\nu_{\s_i}:W_i\to L^2(X_\mathcal{H}))$ of $\mathcal{H}$.
This decomposition induces the spectral decomposition of the
functional $I_\mathcal{F}$.

In fact, the inclusion $X_\mathcal{F}\subset X_\mathcal{H}\subset
X_\mathcal{G}$ induces the period map
$I^{aut}_{\pi,\s_i}\in\mathcal{P}(V,\bc)$ for every $\s_i$, via
the following composition of maps
\begin{equation*}
 I^{aut}_{\pi,\s_i}:V\stackrel{\nu_\pi}{\longrightarrow}C^\8(X_\mathcal{G})
 \stackrel{r_{\mathcal{H}}}{\longrightarrow}C^\8(X_\mathcal{H})
 \stackrel{\nu_{\s_i}^*}{\longrightarrow}W_i
 \stackrel{\nu_{\s_i}}{\longrightarrow}C^\8(X_\mathcal{H})
 \stackrel{r_{\mathcal{F}}}{\longrightarrow}C^\8(X_\mathcal{F})
\stackrel{I_\mathcal{F}}{\longrightarrow}\bc\ .
\end{equation*}
The spectral decomposition \eqref{spec-decomp-H} gives rise to the
decomposition $I_\mathcal{F}=\sum_{\s_i}I^{aut}_{\pi,\s_i}$.
(Since in order to compute $I_\mathcal{F}(v)$ for a vector $v\in
V$, we first can restrict $\nu_\pi(v)$ to the cycle
$X_\mathcal{H}$, decompose it with respect to the action of
$\mathcal{H}$ and then compute the integral over $X_\mathcal{F}$
for each component.) Note that the functional $I^{aut}_{\pi,\s_i}$
is the composition of two automorphic periods:
$I^{aut}_{\pi,\s_i}=T^{aut}_{\s_i,\bc}\circ
T^{aut}_{\pi,\s_i}:V\to W\to\bc$.

We now use the strong Gelfand property for the triple
$\mathcal{F}\subset\mathcal{H}\subset\mathcal{G}$, i.e., that the
product of dimensions $\dim\Hom_\mathcal{H}(V,W_i)\cdot
\dim\Hom_\mathcal{F}(W_i,\bc)\leq 1$ for all $W_i$. We choose
model periods $T^{mod}_{\pi,\s_i}\in \Hom_\mathcal{H}(V,W_i)$,
$T^{mod}_{\s_i,\bc}\in \Hom_\mathcal{F}(W_i,\bc)$. As we explained
above this leads to the automorphic coefficients of
proportionality: $T^{aut}_{\pi,\s_i}=a_{\pi,\s_i}\cdot
T^{mod}_{\pi,\s_i}$ and $T^{aut}_{\s_i,\bc}=b_{\s_i,\bc}\cdot
T^{mod}_{\s_i,\bc}$. We denote by
$\gamma_{\pi,\s_i}=a_{\pi,\s_i}\cdot b_{\s_i,\bc}$ and by
$I^{mod}_{\pi,\s_i}=T^{mod}_{\s_i,\bc}\circ T^{mod}_{\pi,\s_i}\in
V^*$. With such notations we arrive at the spectral decomposition
of the functional $I_\mathcal{F}\in V^*$ which is associated to
the triple of strong Gelfand subgroups
$\mathcal{F}\subset\mathcal{H}\subset\mathcal{G}$
\begin{equation}\label{spec-triple}
 I_\mathcal{F}=\sum_{\s_i\ {\rm automorphic}}
 \gamma_{\pi,\s_i}\cdot I^{mod}_{\pi,\s_i}\ .
\end{equation}

\begin{remark} We note that for a non-compact cycle
$X_\mathcal{F}$ there is no obvious way to write down the analog
of the spectral decomposition \eqref{spec-triple} even if the
initial representation $\pi$ is cuspidal. This is because a priori
the period $T^{aut}_{\s_i,\bc}$ might not be defined for all
$\s_i$ (e.g., non-compact periods of Eisenstein series). Usually,
one has to introduce an appropriate regularization procedure in
order to define the corresponding periods. In this paper we only
consider cycles $X_\mathcal{F}$ which are compact and hence will
not face this problem.
\end{remark}

\subsection{\bf Rankin-Selberg type spectral
identities}\label{RS-subsect} Our main observation is that for a
given pair of groups $\mathcal{F}\subset \mathcal{G}$ there might
be different intermediate subgroups $\mathcal{H}$ as above leading
to {\it different} spectral decompositions of the same functional
$I_\mathcal{F}$ and hence to identities between the automorphic
coefficients.

Let $\mathcal{G}$  be a real Lie group and
$\mathcal{F}\subset\mathcal{H}_i\subset\mathcal{G}$, $i=1,2$, be a
collection of subgroups such that in the following commutative
diagram each embedding is a strong Gelfand pair (i.e., pairs
$(\mathcal{G},\mathcal{H}_i)$ and $(\mathcal{H}_i,\mathcal{F})$
are strong Gelfand pairs)
\begin{equation}\label{comm-1}
\xymatrix{&\mathcal{G}&\\ \mathcal{H}_1\ar@{^{(}->}[ur]& &\mathcal{H}_2
\ar@{_{(}->}[ul]\\
&\mathcal{F}\ar@{^{(}->}[ur]\ar@{_{(}->}[ul]&}
\end{equation}
We call such a collection of subgroups a {\it strong Gelfand
formation}.

Let $\G\subset \mathcal{G}$ be a lattice  and
$X_\mathcal{G}=\G\sm\mathcal{G}$ the corresponding automorphic
space. Let $X_i=X_{\mathcal{H}_i}\subset X_\mathcal{G}$ and
$X_\mathcal{F}\subset X_\mathcal{G}$ be closed orbits of
$\mathcal{H}_i$ and $\mathcal{F}$, respectively, satisfying the
following commutative diagram of embeddings
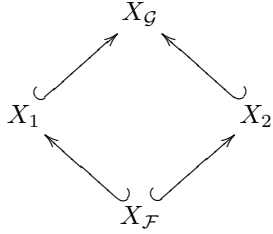
\begin{equation*}\label{comm-2}
\xymatrix{&X_\mathcal{G}&\\ X_1\ar@{^{(}->}[ur]& &X_2\ar@{_{(}->}[ul]\\
&X_\mathcal{F}\ar@{^{(}->}[ur]\ar@{_{(}->}[ul]&}
\end{equation*}
assumed to be compatible with the diagram \eqref{comm-1}. We endow
each orbit (as well as $X_\mathcal{G}$) with a measure invariant
under the corresponding subgroup (to explain the idea, we assume
that all orbits are compact, and hence, these measures could be
normalized to have mass one).

We fix a decompositions $L^2(X_1)=\oplus_i \s_i$ into a direct sum
(in general into a direct integral) of automorphic representations
$(\s_i,\nu_{\s_i},W_i)$ of $\mathcal{H}_1$ and similarly
$L^2(X_2)=\oplus_j \tau_j$ for automorphic representations
$(\tau_j,\nu_{\tau_j}, U_j)$ of $\mathcal{H}_2$.

Let $(\pi,\nu_\pi)$ be an automorphic representation of
$\mathcal{G}$ and $I_\mathcal{F}:V\to\bc$ the period defined by
the integration over the cycle $X_\mathcal{F}$. As we explained in
Section \ref{tripl-Gelf-sect}, two different triples
$\mathcal{F}\subset\mathcal{H}_1\subset\mathcal{G}$ and
$\mathcal{F}\subset\mathcal{H}_2\subset\mathcal{G}$ lead to two
spectral decompositions of the period $I_\mathcal{F}$. Namely, let
us choose the model periods $T^{mod}_{\pi,\s_i}\in
\Hom_{\mathcal{H}_1}(V,W_i)$, $T^{mod}_{\s_i,\bc}\in
\Hom_\mathcal{F}(W_i,\bc)$ and similarly $T^{mod}_{\pi,\tau_j}\in
\Hom_{\mathcal{H}_2}(V,U_j)$, $T^{mod}_{\tau_j,\bc}\in
\Hom_\mathcal{F}(U_j,\bc)$. This leads to the introduction of
constants $T^{aut}_{\pi,\s_i}=a_{\pi,\s_i}T^{mod}_{\pi,\s_i}$,
$T^{aut}_{\s_i,\bc}=b_{\s_i,\bc}T^{mod}_{\s_i,\bc}$ and
$T^{aut}_{\pi,\s_i}=c_{\pi,\tau_j}T^{mod}_{\pi,\tau_j}$,
$T^{aut}_{\tau_j,\bc}=d_{\tau_j,\bc}T^{mod}_{\tau_j,\bc}$ as in
Section \ref{Gelf-period-sect}. We denote by
$\gamma_{\pi,\s_i}=a_{\pi,\s_i}\cdot b_{\s_i,\bc}$,
$I^{mod}_{\pi,\s_i}=T^{mod}_{\s_i,\bc}\circ T^{mod}_{\pi,\s_i}\in
V^*$ and similarly by $\delta_{\pi,\tau_j}=c_{\pi,\tau_j}\cdot
d_{\tau_j,\bc}$, $I^{mod}_{\pi,\tau_j}=T^{mod}_{\tau_j,\bc}\circ
T^{mod}_{\pi,\tau_j}\in V^*$. The spectral decomposition
\eqref{spec-triple} for two triples of Gelfand pairs
\eqref{comm-1} and the corresponding orbits \eqref{comm-2} implies
the following identity
\begin{equation}\label{spec-two-trip}
\sum_{\s_i}\gamma_{\pi,\s_i}\cdot I^{mod}_{\pi,\s_i}=
I_\mathcal{F}=\sum_{\tau_j}\delta_{\pi,\tau_j}\cdot
I^{mod}_{\pi,\tau_j}\ .
\end{equation} We call such an identity the {\it Rankin-Selberg type
spectral identity} or the {\it period identity} associated with
the Gelfand formation $(\mathcal{G},\mathcal{H}_1,
\mathcal{H}_2,\mathcal{F})$, the corresponding orbits and the
automorphic representation $\pi$. Note that the summation on the
left in \eqref{spec-two-trip} is over the set of irreducible
representations of $\mathcal{H}_1$ occurring in $L^2(X_1)$ and the
summation on the right is over the set of irreducible
representations of $\mathcal{H}_2$ occurring in $L^2(X_2)$. Since
groups $\mathcal{H}_1$ and $\mathcal{H}_2$ might be quite
different, the identity \eqref{spec-two-trip} is non-trivial in
general. Surprisingly, even if $\mathcal{H}_1$ and $\mathcal{H}_2$
are conjugate in $\mathcal{G}$, the resulting identity is
non-trivial in general (e.g., for $\mathcal{G}=\PGLR$ and two
unipotent subgroups intersecting over $\mathcal{F}=e$, this gives
the Vorono\u\i\ type summation formula for Fourier coefficients of
cusp forms or of Eisenstein series).

The above identity is the identity between functionals on $V$. It
is easy to translate it into the identities for weighted sums of
coefficients $\gamma$'s and $\delta$'s. Let $v\in V$ be a vector.
It will play the role of a test function. As we explained in
Section \ref{Gelf-period-sect}, in order to construct model
periods, we have to consider model realizations of all
corresponding representations. In particular, we can view $v\in V$
as a function on some manifold (or a section of a vector bundle).
The resulting functionals $I^{mod}_{\pi,\s_i}$ and
$I^{mod}_{\pi,\tau_j}$ could be viewed as integral transforms on
the spaces of such functions. Hence, we obtain for any $v\in V$,
the identity
\begin{equation}\label{spec-two-trip-v}
\sum_{\s_i}\gamma_{\pi,\s_i}\cdot I^{mod}_{\pi,\s_i}(v)=
\sum_{\tau_j}\delta_{\pi,\tau_j}\cdot I^{mod}_{\pi,\tau_j}(v)\
\end{equation} for the weighted sums of products of automorphic
periods $\gamma_{\pi,\s_i}=a_{\pi,\s_i}\cdot b_{\s_i,\bc}$, and
$\delta_{\pi,\tau_j}=c_{\pi,\tau_j}\cdot d_{\tau_j,\bc}$. The main
point of \eqref{spec-two-trip-v} is that the weights
$I^{mod}_{\pi,\s_i}(v)$ and $I^{mod}_{\pi,\tau_j}(v)$ could be
computed in some explicit models without any reference to the
automorphic picture. We will show below that as a special case,
these identities include the classical Rankin-Selberg identity.

\begin{remark}{}We note that one can associate a non-trivial
spectral identity of a kind we described above to a {pair} of
different {\it filtrations} of a group $\mathcal{G}$ by subgroups
forming strong Gelfand pairs. Namely, we can associate a spectral
identity to two filtrations $\{\mathcal{F}= G_0\subset
G_1\subset\dots\subset G_n=\mathcal{G}\}$ and $\{\mathcal{F}=
H_0\subset H_1\subset\dots\subset H_m=\mathcal{G}\}$ of subgroups
in the same group $\mathcal{G}$ such that all pairs
$(G_{i+1},G_i)$ and $(H_{j+1},H_j)$ are strong Gelfand pairs
having the same intersection $\mathcal{F}$. One can also ``twist"
such an identity by a non-trivial character or an automorphic
representation $(\chi, U_\chi)$ of the group $\mathcal{F}$. In
this case the resulting identity is not for an
$\mathcal{F}$-invariant {\it functional} $I_\mathcal{F}$, but for
an automorphic period {\it map} in the period space
$\mathcal{P}(V_\pi,U_\chi)$.
\end{remark}

\subsubsection{Bounds for coefficients}
The Rankin-Selberg type formulas \eqref{spec-two-trip-v} can be
used in order to obtain bounds for the coefficients
$\gamma_{\pi,\s_i}$ and $\delta_{\pi,\tau_j}$ (e.g., Theorems
\ref{thm} and \ref{thm-3}). To this end one has to study
properties of the integral transforms defined by the functionals
$I^{mod}_{\pi,\s_i}$ and $I^{mod}_{\pi,\tau_j}$ on $V$. As
mentioned in Section \ref{Gelf-period-sect}, the construction of
model functionals involves explicit models of representations in
some spaces of functions (or sections of vector bundles). The
model periods $T^{mod}_{\pi,\s_i}$ and $T^{mod}_{\pi,\tau_j}$ are
then given as integral operators with explicit kernels and the
same is true for the resulting model functionals
$I^{mod}_{\pi,\s_i}$ and $I^{mod}_{\pi,\tau_j}$. These functionals
could be defined for all unitary representations $\pi$ of
$\mathcal{G}$, $\s$ of $\mathcal{H}_1$ and $\tau$ of
$\mathcal{H}_2$. Hence we obtain a pair of integral transforms
$h_{\s}=I^{mod}_{\pi,\s}:V^{mod}\to C(\hat{\mathcal{H}}_1)$,
$v\mapsto h_{\s}(v)=I^{mod}_{\pi,\s}(v)$ (here
$\hat{\mathcal{H}}_1$ is the unitary dual of $\mathcal{H}_1$ and
$V^{mod}$ is an explicit model of the representation $V$) and for
the triple $(\mathcal{G},\mathcal{H}_2, \mathcal{F})$ the
transform $g_{\tau}(v)=I^{mod}_{\pi,\tau}(v)$. For the classical
Rankin-Selberg identity this pair of transforms constitutes the
pair of the Fourier and the Mellin transforms on the space of
(smooth with certain decay at infinity) functions on the line
$\br$. The latter is the model for the representation $\pi$ of the
principal series (see Section \ref{unipot-section} for more
details).

For applications, one needs to study analytical properties of
these transforms. This is a problem in harmonic analysis which has
nothing to do with the automorphic picture. We study the
corresponding transforms, in the particular cases under the
consideration, in two technical Lemmas \ref{lem1} and \ref{lem2},
where some instance of what might be called an ``uncertainty
principle" for the pair of such transforms is established.

The idea behind the proofs of Theorems \ref{thm} and \ref{thm-3}
is quite standard by now (and was learned by us from \cite{Go1}).
It is based on the appropriate Rankin-Selberg type identity and
the necessary analytic information for the corresponding integral
transforms (e.g., Lemmas \ref{lem1} and \ref{lem2}). Namely, we
construct a family of test vectors $v_T\in V^{mod} $ parameterizes
by the real parameter $T\geq 1$ such that when substituted into
the Rankin-Selberg type identity \eqref{spec-two-trip-v} it will
pick up the (weighted) sum of coefficients $\gamma_{\pi,\s_i}$ for
$i$ in a certain ``short" interval around $T$ (i.e., the density
$h_\s(v_T)$  is concentrated on $\hat{\mathcal{H}}_1$ around
representations with the parameter of the representation $\s$
close to $T$). We show then that the integral transform
$g_\tau(v_T)$ of such a vector is a slowly changing function on
$\hat{\mathcal{H}}_2$ and estimate its support and the size. This
allows us to bound the right hand side in \eqref{spec-two-trip-v}
using Cauchy-Schwartz inequality and the {\it mean value} (or
convexity) bound for the coefficients $\delta_{\pi,\tau_j}$ (e.g.,
bounds \eqref{convex-in-proof}, \eqref{convex-in-proof-2}). A
simple way to obtain these mean value bounds was explained by us
in \cite{BR3}.

We note that in order to apply this idea to the identity
\eqref{spec-two-trip-v} one needs to have some kind of a {\it
positivity} which is not always easy to achieve. Namely, in order
to bound a single coefficient $\gamma_{\pi,\s_i}$ we have to know
that terms $\gamma_{\pi,\s_i}\cdot I^{mod}_{\pi,\s_i}(v)$ will not
cancel each other in the sum (e.g., all terms are non-negative on
one side of the identity). As a result of this constraint there
are many identities from which it is not clear how to deduce
bounds for the corresponding coefficients. In the examples that we
consider in this paper we choose representations
$V=\mathcal{V}\otimes\mathcal{\bar V}$ of the group
$\mathcal{G}=G\times G$ with $\mathcal{V}$ an irreducible unitary
representation of some other group $G$. For such representations
the necessary positivity is automatic.

In this paper we implement the above strategy for $G=\PGLR$ and
two cases: for the unipotent subgroup $N\subset G$ and for a
compact subgroup $K\subset G$. The first case corresponds to the
unipotent Fourier coefficients and the formula we obtain is
equivalent to the classical Rankin-Selberg formula. The second
case corresponds to the spherical Fourier coefficients which were
introduced by H. Peterson a long time ago, but the corresponding
formula (see Theorem \ref{thm-2}) has never appeared in print, to
the best of our knowledge.

To relate these cases to the above discussion of Rankin-Selberg
type spectral formulas, we set $\mathcal{G}=G\times G$,
$\mathcal{H}_2=\Dl G\ {\hookrightarrow}\ G\times G$ in both cases
under the consideration and $\mathcal{H}_1=N\times N$,
$\mathcal{F}=\Dl N\ {\hookrightarrow}\ N\times N\
{\hookrightarrow}\ G\times G$ for the unipotent Fourier
coefficients and $\mathcal{H}_1=K\times K$, $\mathcal{F}=\Dl K
\hookrightarrow K\times K \hookrightarrow G\times G$ for the
spherical Fourier coefficients. Strictly speaking,  the uniqueness
principle is only ``almost'' satisfied for the subgroup $N$, but
the theory of the constant term of the Eisenstein series provides
the necessary remedy in the automorphic setting (see Section
\ref{triple}).

We also illustrate {\it analytic applications} of these
Rankin-Selberg type spectral identities. We prove non-trivial
bounds for both types of these Fourier coefficients. While bounds
for the unipotent coefficients (Theorem \ref{thm}) are known (and
even much better bounds are known for the Hecke-Maass forms), for
the spherical case our bounds (Theorem \ref{thm-3}) are new. As a
corollary, we obtain a subconvexity bound for certain automorphic
$L$-functions.

The method described above also lies behind the subconvexity for
the triple $L$-function given in \cite{BR4}. There the
corresponding strong Gelfand formation consists of
$\mathcal{G}=G\times G\times G\times G$ with $G=\PGLR$,
$\mathcal{F}=\Dl G$ and $H_i= G\times G$ with two different
embedding into $\mathcal{G}$. Recently it became clear that there
are many strong Gelfand formations in higher rank groups. We hope
to discuss the corresponding identities elsewhere.

\subsection{Unipotent Fourier coefficients of Maass forms}

Let $G=\PGLR$ and denote by $K=PO(2)$ the standard maximal compact
subgroup of $G$. Let $\uH=G/K$ be the upper half plane endowed
with a hyperbolic metric and the corresponding volume element
$d\mu_\uH$.

Let $\G\subset G$ be a non-uniform lattice.  We assume for
simplicity that, up to equivalence, $\G$ has a unique cusp which
is reduced at $\8$. This means that the unique up to conjugation
unipotent subgroup
$\G_\8\subset \G$ is generated by $\left(\begin{array}{cc} 1&  1\\
0& 1\end{array}\right)$ (e.g., $\G=PSL_2(\bz)$). We denote by
$X=\G\sm G$ the automorphic space and by $Y=X/K=\G\sm\uH$ the
corresponding Riemann surface (with possible conic singularities
if $\G$ has elliptic elements). This induces the corresponding
Riemannian metric on $Y$, the volume element $d\mu_Y$ and the
Laplace-Beltrami operator $\Dl$. We normalize $d\mu_Y$ to have the
total {\it volume one}.

Let $\phi_\tau\in L^2(Y)$ be a Maass cusp form. In particular,
$\phi_\tau$ is an eigenfunction of $\Dl$ with the eigenvalue which
we write in the form $\mu=\frac{1-\tau^2}{4}$ for some $\tau\in
\bc$. We will always assume that $\phi_\tau$ is normalized to have
{\it $L^2$-norm one}. We can view $\phi_\tau$ as a $\G$-invariant
eigenfunction of the Laplace-Beltrami operator $\Dl$ on $\uH$.
Consider the classical Fourier expansion of $\phi_\tau$ at $\8$
given by (see \cite{Iw})
\begin{equation}\label{f-c}
 \phi_\tau(x+iy)=\sum_{n\neq
0}a_n(\phi_\tau)\CW_{\tau,n}(y)e^{2\pi in x}\ .
\end{equation}
Here $\CW_{\tau,n}(y)e^{2\pi in x}$ are properly normalized
eigenfunctions of $\Dl$ on $\uH$ with the same eigenvalue $\mu$ as
that of the function $\phi_\tau$. The functions $\CW_{\tau,n}$ are
usually described in terms of the $K$-Bessel function. In Section
\ref{u-fourier} we recall the well-known description of functions
$\CW_{\tau,n}$ in terms of certain matrix coefficients for unitary
representations of $G$.

We note that from the group-theoretic point of view, the Fourier
expansion \eqref{f-c} is a consequence of the decomposition of the
function $\phi_\tau$ under the natural action of the group
$N/\G_\8$ (commuting with $\Dl$). Here $N$ is the standard
upper-triangular subgroup and the decomposition is with respect to
the characters of the group $N/\G_\8$ (see Section
\ref{u-fourier}).

The vanishing of the zero Fourier coefficient $a_0(\phi_\tau)$ in
\eqref{f-c} distinguishes cuspidal Maass forms (for $\G$ having
several inequivalent cusps, the vanishing of the zero Fourier
coefficient is required at each cusp).

The coefficients $a_n(\phi_\tau)$ are called the Fourier
coefficients of the Maass form $\phi_\tau$ and play a prominent
role in analytic number theory.

One of the central problems in the analytic theory of automorphic
functions is the following

{\bf Problem:} Find the best possible constants $\s$, $\rho$ and
$C_\G$ such that the following bound holds
$$|a_n(\phi_\tau)|\leq C_\G\cdot |n|^\s\cdot(1+|\tau|)^\rho\ .$$
In particular, one asks for constants $\s$ and $\rho$ which are
{\it independent} of $\phi_\tau$ (i.e., depend on $\G$ only; for a
brief discussion of the question history, see Remark \ref{his}).

It is easy to obtain a polynomial bound for coefficients
$a_n(\phi_\tau)$ using boundedness of $\phi_\tau$ on $Y$. Namely,
G. Hardy and E. Hecke essentially proved that the following bound
$$\sum_{|n|\leq
T}|a_n(\phi_\tau)|^2\leq C\cdot\max\{ T, 1+|\tau|\}$$ holds for
any $ T\geq 1$,  with the constant $C$ depending on $\G$ only (see
\cite{Iw}). It would be very interesting to improve this bound for
coefficients $a_n(\phi_\tau)$ in the range $|n|\ll|\tau|$.

For a fixed $\tau$, we have the bound $|a_n(\phi_\tau)|\leq
C_\tau|n|^\haf$. This bound is usually called the standard bound
or the Hardy/Hecke bound for the Fourier coefficients of cusp
forms (in the $n$ aspect).

The first improvements of the standard bound are due to H.
Sali{$\rm\acute{e}$} and A. Walfisz using exponential sums. Rankin
\cite{Ra} and Selberg \cite{Se} independently discovered the
so-called Rankin-Selberg unfolding method (i.e., the formula
\eqref{unfold} below) which allowed them to show that for any
$\eps>0$, the bound $|a_n(\phi)|\ll|n|^{\frac{3}{10}+\eps}$ holds.
Their approach is based on the integral representation for the
weighted sum of Fourier coefficients $a_n(\phi)$. To state it, we
assume, for simplicity, that the so-called residual spectrum is
trivial (i.e., the Eisenstein series $E(s,z)$ are holomorphic for
$s\in(0,1)$; e.g, $\G=PGL_2(\bz)$). (The reader also should keep
in mind that we use the normalization ${\rm vol}(Y)=1$ and ${\rm
vol}(\G_\8\sm N)=1$.) We have then
\begin{eqnarray}\label{RS-Fourier}
\sum_n|a_n(\phi)|^2\hat\al(n)=\al(0)+\frac{1}{2\pi
i}\int_{Re(s)=\haf} D(s,\phi,\bar\phi)M(\al)(s)ds\ ,
\end{eqnarray}
where $\al\in C^\8(\br)$ is an appropriate test function with the
Fourier transform $\hat\al$ and the Mellin transform $M(\al)(s)\
$,
\begin{eqnarray}\label{rs}
D(s,\phi,\bar\phi)=\G(s,\tau)\
\cdot\langle\phi\bar\phi,E(s)\rangle_{L^2(Y)}\ ,
\end{eqnarray} where
$E(z,s)$ is an appropriate non-holomorphic Eisenstein series and
$\G(s,\tau)$ is given explicitly  in terms of the Euler
$\G$-function (see Remark \ref{his}).

The proof of \eqref{RS-Fourier}, given by Rankin and Selberg, is
based on the so-called unfolding trick, which amounts to the
following. Let $E(s,z)$ be the Eisenstein series given by
$E(s,z)=\sum_{\g\in\G_\8\sm\G}y^s(\g z)$ for $Re(s)>1$ (and
analytically continued to a meromorphic function for all
$s\in\bc$). We have the following ``unfolding" identity valid for
$Re(s)>1$,
\begin{equation}
\begin{split}\label{unfold}
&\langle\phi\bar\phi,E(z,s)\rangle_{L^2(Y)}=
\int_{\G\sm\uH}\phi(z)\bar\phi(z)\sum_{\g\in\G_\8\sm\G}y^s(\g
z)d\mu_Y=\\
&=\int_{\G_\8\sm\uH}\phi(z)\bar\phi(z)y^s(z)d\mu_\uH=
\int_0^\8\left(\int_0^1\phi(x+iy)\bar\phi(x+iy)\ dx\right)
y^{s-1}\ d^xy\ .
\end{split}\end{equation}
The Mellin inversion formula, together with the Fourier expansion
\eqref{f-c} for $\phi$, leads to the Rankin-Selberg formula
\eqref{RS-Fourier}.

Using the strategy formulated in Section \ref{RS-subsect}, in this
paper we deduce the Rankin-Selberg formula \eqref{RS-Fourier}
directly from the uniqueness principle in representation theory
and hence avoid the use of the unfolding trick \eqref{unfold}. One
of the uniqueness results we are going to use is related to the
unipotent subgroup $N\subset G$ such that $\G_\8\subset N$ (the
so-called $\G$-cuspidal unipotent subgroup). In fact, the
definition of classical Fourier coefficients $a_n(\phi_\tau)$ is
implicitly based on the uniqueness of $N$-equivariant functionals
on an irreducible (admissible) representation of $G$ (i.e., on the
uniqueness of the so-called Whittaker functional). For this
reason, we call the coefficients $a_n(\phi_\tau)$ the {\it
unipotent} Fourier coefficients.

We obtain a somewhat different (a slightly more ``geometric") form
of the Rankin-Selberg identity \eqref{RS-Fourier}. In particular,
we exhibit a connection between analytic properties of the
function $D(s,\phi,\bar\phi)$ and analytic properties of certain
invariant functionals on irreducible unitary representations of
$G$. This allows us to deduce subconvexity bounds for Fourier
coefficients of Maass forms for a {\it general} $\G$ in a more
transparent way (here we relay on ideas of A. Good \cite{Go1} and
on our earlier results \cite{BR1} and \cite{BR3}). Namely, we
prove the following bound for the Fourier coefficients
$a_n(\phi_\tau)$.

\begin{theorem}\label{thm}
Let $\phi_\tau$ be a fixed Maass form of $L^2$-norm one. For any
$\eps>0$, there exists an explicit constant $C_\eps$ such that
$$\sum\limits_{|k-T|\leq T^{\frac{2}{3}}}|a_k(\phi_\tau)|^2\leq
C_\eps\cdot T^{\frac{2}{3}+\eps}\ .$$
\end{theorem}

In particular, we have $|a_n(\phi_\tau)|\ll
|n|^{\frac{1}{3}+\eps}$. This is weaker than the Rankin-Selberg
bound, but holds for {\it general lattices} $\G$ (i.e., not
necessary a congruence subgroup). The bound in the theorem was
first claimed in \cite{BR1} and the analogous bound for
holomorphic cusp forms was proved by Good \cite{Go1} by a
different method. Here we give full details of the proof following
a slightly different argument.

The main goal of this paper, however, is different. Our main new
results deal with another type of Fourier coefficients associated
with a Maass form. These Fourier coefficients, which we call
spherical, are associated to a compact subgroup of $G$.

\subsection{Spherical Fourier coefficients}\label{anisot-section}
When dealing with spherical Fourier coefficients we assume, for
simplicity, that $\G\subset G$ is a co-compact subgroup and
$Y=\G\sm \uH$ is the corresponding compact Riemann surface. Let
$\phi_\tau$ be a norm one eigenfunction of the Laplace-Beltrami
operator on $Y$, i.e., a Maass form. We would like to consider a
kind of a Taylor series expansion for $\phi_\tau$ at a point on
$Y$. To define this expansion, we view $\phi_\tau$ as a
$\G$-invariant eigenfunction on $\uH$. We fix a point $z_0\in
\uH$. Let $z=(r,\theta)$, $r\in\br^+$ and $\theta\in S^1$, be the
geodesic polar coordinates centered at $z_0$ (see \cite{He}). We
have the following spherical Fourier expansion of $\phi_\tau$
associated to the point $z_0$
\begin{equation}\label{spherical-expansion}
 \phi_\tau(z)=\sum_{n\in\bz}b_{n,z_0}(\phi_\tau)P_{\tau,n}(r)e^{ in \theta}\
 .
\end{equation}
Here functions $P_{\tau,n}(r)e^{ in \theta}$ are properly
normalized eigenfunctions of $\Dl$ on $\uH$ with the same
eigenvalue $\mu$ as that of the function $\phi_\tau$. The
functions $P_{\tau,n}$ can be described in terms of the classical
Gauss hypergeometric function or the Legendre function. In Section
\ref{P-n}, we will describe special functions $P_{\tau,n}$ and
their normalization in terms of certain matrix coefficients of
irreducible unitary representations of $G$.

We call the coefficients $b_n(\phi_\tau)=b_{n,z_0}(\phi_\tau)$ the
spherical (or anisotropic) Fourier coefficients of $\phi_\tau$
(associated to a point $z_0$). These coefficients were introduced
by H.~Petersson and played a major role in recent works of Sarnak
(e.g., \cite{Sa3}). It was discovered by J.-L. Waldspurger
\cite{Wa} that in certain cases these coefficients are related to
special values of $L$-functions (see Remark \ref{L-funct}).

As in the case of the unipotent expansion \eqref{f-c}, the
spherical expansion \eqref{spherical-expansion} is the result of
an expansion with respect to a  group action. Namely, the
expansion \eqref{spherical-expansion} is with respect to
characters of the compact subgroup $K_{z_0}={\rm Stab}_{z_0}G$
induced by the natural action of $G$ on $\uH$ (for more details,
see Section \ref{sperical-section}).

The expansion \eqref{spherical-expansion} exists for {\it any}
eigenfunction of $\Dl$ on $\uH$ . This follows from a simple
separation of variables argument applied to the operator $\Dl$ on
$\uH$. For a proof and a discussion of the growth properties of
coefficients $b_n(\phi)$ for a general eigenfunction $\phi$ on
$\uH$, see \cite{He}, \cite{L}. For another approach which is
applicable to Maass forms, see \cite{BR2}.

Under the normalization we choose, the coefficients
$b_n(\phi_\tau)$ are bounded on the average. Namely, one can show
that the following bound holds
$$\sum_{|n|\leq T}|b_n(\phi_\tau)|^2\leq  C'\cdot\max\{ T,
1+|\tau|\}\ $$ for any $T\geq 1$, with the constant $C'$ depending
on $\G$ only (see \cite{R}).

As our approach is based directly on the uniqueness principle, we
are able to prove an analog of the Rankin-Selberg formula
\eqref{RS-Fourier} with the group $N$ replaced by a maximal
compact subgroup of $G$. This is the main aim of the paper. We
obtain an analog of the Rankin-Selberg formula \eqref{RS-Fourier}
for the coefficients $b_n(\phi_\tau)$. Roughly speaking, new
formula amounts to the following (for the exact form, see formula
\eqref{rs-K-u})

\begin{theorem}\label{thm-2} Let $\{\phi_{\lm_i}\}$ be an
orthonormal basis of $L^2(Y)$ consisting of Maass forms. Let
$\phi_\tau$ be a fixed Maass form.

There exists an explicit integral transform $^\sharp:C^\8(S^1)\to
C^\8(\bc)$, $u(\theta)\mapsto u^\sharp_\tau(\lm)$, such that for
all $u\in C^\8(S^1)$, the following relation holds:
\begin{eqnarray}\label{rs-K-u-crude}
\sum_n|b_n(\phi_\tau)|^2 \hat u(n)= u(1)+
\sum_{\lm_i\not=1}\mathcal{L}_{z_0}(\phi_{\lm_i})\cdot
u^\sharp_\tau(\lm_i)\ ,
\end{eqnarray}
with some explicit coefficients
$\mathcal{L}_{z_0}(\phi_{\lm_i})\in \bc$ which are independent of
$u$.

Here $\hat u(n)=\frac{1}{2\pi}\int\limits_{S^1}
u(\theta)e^{-in\theta}d\theta$ and $u(1)$ is the value at $1\in
S^1$.
\end{theorem}

The definition of the integral transform $^\sharp$ is based on the
uniqueness of certain invariant trilinear functionals on
irreducible unitary representations of $G$ and described
explicitly in the formula \eqref{sperical-transf}. The trilinear
functional was studied by us in \cite{BR3} and \cite{BR4}. The
main point of the relation \eqref{rs-K-u-crude} is that the
transform $u^\sharp_\tau(\lm_i)$ depends only on the {\it
parameters} $\lm_i$ and $\tau$, but {\it not} on the choice of
Maass forms $\phi_{\lm_i}$ and $\phi_\tau$. The coefficients
${\CL}_{z_0}(\phi_{\lm_i})$ are essentially given by the product
of the triple product coefficients
$\langle\phi_\tau^2,\phi_{\lm_i}\rangle_{L^2(Y)}$ and the values
of Maass forms $\phi_{\lm_i}$ at the point $z_0$. In some special
cases both types of these coefficients are related to
$L$-functions (see \cite{W}, \cite{JN}, \cite{Wa} and Remark
\ref{L-funct}).

A formula similar to \eqref{rs-K-u-crude} holds for a non-uniform
lattice $\G$ as well, and includes the contribution from the
Eisenstein series (see the formula in Remark \ref{rem-sph}). Also,
a similar formula holds for holomorphic forms. We intend to
discuss it elsewhere.

The new formula \eqref{rs-K-u-crude} allows us to deduce the
following bound for the {\it spherical} Fourier coefficients of
Maass forms.

\begin{theorem}\label{thm-3}
Let $\G$ be as above and $\phi_\tau$ is a fixed Maass form of
$L^2$-norm one. For any $\eps>0$, there exists an explicit
constant $D_\eps$ such that $$\sum_{|k-T|\leq
T^{\frac{2}{3}}}|b_k(\phi_\tau)|^2\leq D_\eps\cdot
T^{\frac{2}{3}+\eps}\ .$$
\end{theorem}

In particular, we have $|b_n(\phi_\tau)|\ll
|n|^{\frac{1}{3}+\eps}$ for any $\eps>0$. An analogous bound
should hold for spherical Fourier coefficients of holomorphic cusp
forms. We hope to return to this subject elsewhere.

The proof of the bound in the theorem follows from essentially the
same argument as in the case of the unipotent Fourier
coefficients, once we have the Rankin-Selberg type identity
\eqref{rs-K-u-crude}. In the proof we use bounds for triple
products of Maass forms obtained in \cite{BR3}, and a well-known
bound for the averaged value of eigenfunctions of $\Dl$.

In special cases, the bound in the theorem could be interpreted as
a subconvexity bound for some automorphic $L$-function (see Remark
\ref{L-funct}).

\subsection{Remarks}
\subsubsection{Special values of $L$-functions}\label{L-funct}
One of the reasons one might be interested in bounds for
coefficients $b_k(\phi_\tau)$ is their relation to certain
automorphic $L$-functions. It was discovered by J.-L. Waldspurger
\cite{Wa} that, in certain cases, the coefficients
$b_k(\phi_\tau)$ are related to special values of $L$-functions.
H. Jacquet constructed the appropriate relative trace formula
which covers these cases (see \cite{JN}). The simplest case of the
formula of Waldspurger is the following. Let $z_0= i\in
SL_2(\bz)\sm\uH$ and $E=\mathbb{Q}(i)$. Let $\pi$ be the
automorphic representation which corresponds to $\phi_\tau$, $\Pi$
its base change over $E$ and $\chi_n(z)=\left(z/\bar
z\right)^{4n}$ the $n$-th power of the basic Gr\"ossencharacter of
$E$. One has then, under appropriate normalization (for details,
see \cite{Wa}, \cite{JN}, \cite{MW}), the following beautiful
formula
\begin{eqnarray}\label{Wald}
|b_n(\phi_\tau)|^2=\frac{L(\haf,\Pi\otimes\chi_n)}{L(1,Ad\pi)}\ \
.
\end{eqnarray} Using this formula, we can interpret the bound in
Theorem \ref{thm-3} as a bound on the corresponding $L$-functions.
In particular, we have $|L(\haf,\Pi\otimes\chi_n)|\ll
|n|^{2/3+\eps}$. This gives a subconvexity bound (with the
convexity bound for this $L$-function being
$|L(\haf,\Pi\otimes\chi_n)|\ll|n|^{1+\eps}$). The exponent in the
bound corresponds to what is known as a H. Weyl type subconvexity
bound for an $L$-function.

The subconvexity problem is a classical question in analytic
theory of $L$-functions which received much of attention in recent
years (we refer to the survey \cite{IS} for the discussion of
subconvexity for automorphic $L$-functions). In fact, Y. Petridis
and P. Sarnak \cite{PS} recently considered more general
$L$-functions. Among other things, they have shown that
$|L(\haf+it_0,\Pi\otimes\chi_n)|\ll |n|^{\frac{159}{166}+\eps}$
for any fixed $t_0\in\br$ and any automorphic cuspidal
representation $\Pi$ of $GL_2(E)$ (not necessarily a base change).
Their method is also spectral in nature although it uses
Poincar\'{e} series and treats $L$-functions through (unipotent)
Fourier coefficients of cusp forms. We deal directly with periods
and the special value of $L$-functions only appear through the
Waldspurger formula. Of course, our interest in Theorem
\ref{thm-3} lies not so much in the slight improvement of the
Petridis-Sarnak bound for these $L$-functions, but in the fact
that we can give a general bound valid for {\it any} point $z_0$.
(It is clear that for a generic point or a cusp form which is not
a Hecke form, coefficients $b_n$ are not related to special values
of $L$-functions.)

Recently, A. Venkatesh \cite{V} announced (among other remarkable
results) a slightly weaker subconvexity bound for coefficients
$b_n(\phi_\tau)$ for a fixed $\phi_\tau$. His method seems to be
quite different and is based on ergodic theory. In particular, it
is not clear how to deduce the identity \eqref{rs-K-u-crude} from
his considerations. On the other hand, the ergodic method gives a
bound for Fourier coefficients for higher rank groups (e.g., on
$GL(n)$) while it is not yet clear in what higher-rank cases one
can develop Rankin-Selberg type formulas similar to
\eqref{rs-K-u-crude} which would lead to bounds for the
corresponding coefficients.

\subsubsection{Fourier expansions along closed geodesics} There is
one more case where we can apply the uniqueness principle to a
subgroup of $\PGLR$. Namely, we can consider closed orbits of the
diagonal subgroup $A\subset\PGLR$ acting on $X$. It is well-known
that such an orbit corresponds to a closed geodesic on $Y$ (or to
a geodesic ray starting and ending at cusps of $Y$). Such closed
geodesics give rise to Rankin-Selberg type formulas similar to
ones we considered for closed orbits of subgroups $N$ and $K$. In
special cases the corresponding Fourier coefficients are related
to special values of various $L$-functions (e.g., the standard
Hecke $L$-function of a Hecke-Maass forms which appears for a
geodesic connecting cusps of a congruence subgroup of
$PSL(2,\bz)$). In fact, in the language of representations of
ad\`ele groups, which is appropriate for arithmetic $\G$, the case
of closed geodesics corresponds to real quadratic extensions of
$\mathbb{Q}$ (e.g., twisted periods along Heegner cycles) while
the anisotropic expansions (at Heegner points) which we considered
in Section \ref{anisot-section} correspond to imaginary quadratic
extensions of $\mathbb{Q}$ (e.g., twisted ``periods" at Heegner
points).

In order to prove an analog of Theorems \ref{thm} and \ref{thm-3}
for the  Fourier coefficients associated to a closed geodesic, one
has to face certain technical complications. Namely, for orbits of
the diagonal subgroup $A$ one has to consider contributions from
representations of discrete series, while for subgroups $N$ and
$K$ this contribution vanishes. It is more cumbersome to compute a
contribution from discrete series as these representations do not
have nice geometric models. Hence, while the proof of an analog of
Theorem \ref{thm-2} for closed geodesics is straightforward, one
has to study invariant trilinear functionals on discrete series
representations more closely in order to deduce bounds for the
corresponding coefficients. We hope to return to this subject
elsewhere.
\subsubsection{Dependence on the eigenvalue} From the proof that we
present it follows that the constants $C_\eps$ and $D_\eps$ in
Theorems \ref{thm-2} and \ref{thm-3} satisfy the following bound
$$C_\eps,\ D_\eps\leq C(\G)\cdot(1+|\tau|)\cdot|\ln\eps|\ ,$$
for any $0<\eps\leq 0.1$, and some explicit constant $C(\G)$
depending on the lattice $\G$ only.
\subsubsection{Historical remarks}\label{his} The question of the size of
Fourier coefficients of cusp forms was posed (in the $n$ aspect)
by S. Ramanujan for holomorphic forms (i.e., the celebrated
Ramanujan conjecture established in full generality by P. Deligne
for the holomorphic Hecke cusp form for congruence subgroups) and
was extended by H. Petersson to include Maass forms (i.e., the
Ramanujan-Petersson conjecture for Maass forms). In recent years
the $\tau$ aspect of this problem also turned out to be important.

Under the normalization that we have chosen, it is expected that
the coefficients $a_n(\phi_\tau)$ are at most slowly growing as
$n\to\8$ (\cite{Sa3}). Moreover, it is quite possible that the
strong uniform bound of the form $|a_n(\phi_\tau)|\ll
(|n|(1+|\tau|))^\eps$ holds for any $\eps>0$ (e.g.,
Ramanujan-Petersson conjecture for Hecke-Maass forms for
congruence subgroups of $PSL_2(\bz)$). We note, however, that the
behavior of Maass forms and holomorphic forms in these questions
might be quite different (e.g., high multiplicities of holomorphic
forms).

Using the integral representation \eqref{RS-Fourier} and detailed
information about Eisenstein series available only for {\it
congruence subgroups}, Rankin and Selberg showed that for a cusp
form $\phi$ for a congruence subgroup of $PGL_2(\bz)$ one has
$\sum_{|n|\leq T}|a_n(\phi)|^2= CT+O(T^{3/5+\eps})$ for any
$\eps>0$. In particular, this implies that for any $\eps>0$,
$|a_n(\phi)|\ll|n|^{\frac{3}{10}+\eps}$. Since their
groundbreaking papers, this bound was improved many times by
various methods (with the current record for Hecke-Maass forms
being $7/64\approx 0.109...$ due to H. Kim, F. Shahidi and P.
Sarnak \cite{KSa}).

The approach of Rankin and Selberg  is based on the integral
representation of the Dirichlet series given for $Re(s)>1$, by the
series $D(s,\phi,\bar\phi)=\sum_{n>0}\frac{|a_n(\phi)|^2}{n^s}\ .$

The introduction of the so-called Ranking-Selberg $L$-function
$L(s,\phi\otimes\bar\phi)=\zeta(2s)D(s,\phi,\bar\phi)$ played an
even more important role in the further development of automorphic
forms than the bound for Fourier coefficients which Rankin and
Selberg obtained.

Using  integral representation \eqref{rs}, Rankin and Selberg
analytically continued the function $L(s,\phi\otimes\bar\phi)$ to
the whole complex plane and obtained an effective bound for the
function $L(s,\phi\otimes\bar\phi)$ on the critical line
$s=\haf+it$ for $\G$ being a congruence subgroup of $SL_2(\bz)$.
From this, using standard methods in the theory of Dirichlet
series, they were able to deduce  bounds for Fourier coefficients
of cusp forms. In fact, Rankin and Selberg appealed to the
classical Perron formula (in the form given by E. Landau) which
relates analytic behavior of a Dirichlet series with non-negative
coefficients to partial sums of its coefficients. The necessary
analytic properties of $L(s,\phi\otimes\bar\phi)$ are inferred
from properties of the Eisenstein series through the formula
\eqref{rs}.

A small drawback of the original Rankin-Selberg argument is that
their method is applicable to Maass (or holomorphic) forms coming
from {\it congruence subgroups} only. The reason for such a
restriction is the absence of methods which would allow one to
estimate unitary Eisenstein series for general lattices $\G$.
Namely, in order to effectively use the Rankin-Selberg formula
\eqref{RS-Fourier} one would have to obtain {\it polynomial}
bounds for the normalized inner product
$D(s,\phi,\bar\phi)=\G(s,\tau)\cdot\langle\phi\bar\phi,E(s)\rangle_{L^2(Y)}$.
This turns out to be notoriously difficult because of the {\it
exponential} growth of the factor
$\G(s,\tau)=\frac{2\pi^s\G(s)}{\G^2(s/2)\G(s/2+\tau/2)\G(s/2-\tau/2)}$,
for $|s|\to\8$, $s\in i\br$. For a congruence subgroup, the
question could be reduced to known bounds for the Riemann zeta
function or for Dirichlet $L$-functions, as was shown by Rankin
and Selberg. The problem of how to treat general $\G$ was posed by
Selberg in his celebrated paper \cite{Se}.

The breakthrough in this direction was achieved in works of Good
\cite{Go1} (for holomorphic forms) and Sarnak \cite{Sa3} (in
general) who proved non-trivial bounds for Fourier coefficients of
cusp forms for a general $\G$ using spectral methods. The method
of Sarnak was finessed in \cite{BR1} by introducing various ideas
from the representation theory and further extended in \cite{KS}.
The method of this paper is different and avoids the use of
analytic continuation which is central for \cite{Sa3}, \cite{BR1}
and \cite{KS}.

Special cases of the Rankin-Selberg spectral identities described
in the introduction were obtained before by a different method.
The first (vague) attempt to write the above-mentioned formula for
four copies of $G=\PGLR$ and all representations coming from
Eisenstein series was made by N. Kuznetsov \cite{Kz}. His aim was
to obtain a formula representing the eighth moment of the Riemann
zeta function. Later Y. Motohashi \cite{Mo1} obtained the formula
for the fourth moment of the Riemann zeta function. This
corresponds to our identity with $\mathcal{G}=G\times G$ for
$G=\PGLR$, $H_1=T\times T$ where $T\subset G$ is the diagonal
subgroup, $\mathcal{F}=\Dl T$ and $H_2=\Dl G$. To obtain the
fourth moment of the Riemann zeta function, Motohashi considers
representations coming from Eisenstein series. This leads to
considerable technical difficulties which one should not
underestimate. Both Kuznetsov and Motohashi based their approach
on the celebrated Bruggeman-Kuznetsov trace formula (applying it
twice!). The setup we present here even if it does not simplify
the arguments at least gives a more conceptual explanation for the
terms appearing in these identities.

Many other cases of these identities appeared more recently
(mostly stated implicitly as a tool for estimation of
$L$-functions or other quantities). Among these are works of R.
Bruggeman, V. Bykovski\u\i, A.  Ivi\'c,  M. Jutila, P. Michel, A.
I. Vinogradov to name a few.

Finally, we would like to mention that recently R. Bruggeman, M.
Jutila and Y. Motohashi (see \cite{BM}, \cite{Mo2} and references
therein) developed what they call the inner product method. It is
based on the unfolding of an appropriate Poincar\'{e} or Petersson
type series. The standard unfolding leads to the spectral
expansion for the series of the type $\sum_k
A_k(\phi)A_{k+h}(\bar\phi)W(k)$, where $\phi$ is a Maass form and
$A_k$ are appropriate Fourier coefficients (e.g., unipotent or
spherical Fourier coefficients we discussed above). The formulas
obtained in such a way are special cases of our Rankin-Selberg
type formula \eqref{RS4} for a {\it special} test vectors $v$.
These vectors are constructed from certain functions on the
upper-half plane. As a result, the corresponding weights in the
Rankin-Selberg type formulas are reminiscent of {\it exponential}
weights considered by Selberg and Rankin. It seems that using our
approach one can avoid a difficult task of removing these unwanted
weights.

%%%%%%%%%%%%%%%%%%%%%%%%%%%%%%%%%%%%%%%%

The paper is organized as follows. In Section \ref{reps-sect}, we
quickly recall the notion of automorphic representations of $G$
and describe the standard models of representations we will use.

In Section \ref{unipot-section} we reprove the classical
Rankin-Selberg formula and deduce bounds for the unipotent Fourier
coefficients of Maass forms. The proof is based on the uniqueness
of trilinear invariant functionals on irreducible unitary
representations of $G$. We use the description of these
functionals obtained in \cite{BR3}.

In Section \ref{sperical-section} we apply the same strategy to
the spherical Fourier coefficients. In fact, in this case the
proof is less involved since we do not need the theory of the
Eisenstein series in order to remedy the non-uniqueness of
$N$-invariant functionals on irreducible representations of $G$.
Section \ref{sperical-section} contains our main new results and
the reader might read this section independently of Section
\ref{unipot-section}.

In the appendix we prove an asymptotic expansion of the model
trilinear functional. We use this analysis in the proof of Theorem
\ref{thm-3}.

%%%%%%%%%%%%%%%%%%%%%%%%%%%%%%%%%%%%%%%%%%%%%%%%%
{\bf Acknowledgments.} This paper is a byproduct of a joint work
with Joseph Bernstein. It was written under his insistence. It is
a special pleasure to thank him for numerous discussions and for
his constant encouragement and support over many years. I also
would like to thank Peter Sarnak for stimulating discussions and
support.

\section{Representations of $\PGLR$}\label{reps-sect}

We start with a short reminder about the connection between Maass
forms and representation theory of $\PGLR$ which is due to Gelfand
and Fomin (see \cite{GGPS}).

\subsection{Models of representations} All irreducible unitary
representations of the group $G=\PGLR$ are classified. For
simplicity we consider those with a nonzero $K$-fixed vector
(so-called representations of class one) since only these
representations arise from Maass forms. These are the
representations of the principal and the complementary series and
the trivial representation. We will use the following standard
model (or realization) for these representations.

For every complex number $\tau$, consider the space $V_\tau$ of
smooth even homogeneous functions on $\br^2\sm0$ of the
homogeneous degree $\tau-1 \ $ (which means that
$f(ax,ay)=|a|^{\tau-1}f(x,y)$ for all $a\in\br \setminus 0$). The
representation $(\pi_\tau, V_\tau)$ is induced by the action of
the group $\GLR$ given by the formula $\pi_\tau (g)f(x,y)= f(g\inv
(x,y))|\det g|^{\frac{\tau-1}{2}}$. This action is trivial on the
center of $\GLR$ and hence defines a representation of $G$. The
representation $(\pi_\tau, V_\tau)$ is called {\it representation
of the generalized principal series}.

For explicit computations it is often convenient to pass from the
plane model to a line model. Namely, the  restriction of functions
in $V_\tau$ to the line $(x,1) \subset \br^2$ defines an
isomorphism of the space $V_\tau$ with the space $C^\8_\tau(\br)$
of restrictions of smooth homogeneous functions (e.g., decaying at
infinity as $|x|^{\tau-1}$). Hence we can think about vectors in
$V_\tau$ as functions on $\br$.

In the line model the action of an element
$\tilde{a}=diag(a,a\inv)$, $a\in\br^\times$, in the diagonal
subgroup is given by $
\pi_\tau\left(\tilde{a}\right)f(x,1)=f(a\inv x,
a)=|a|^{\tau-1}f(a^{-2}x,1)\ ;$
and the action of an element $\tilde{n}=\left(\begin{array}{cc}1& n \\
& 1\end{array}\right)$ in the unipotent group is given  $
\pi_\tau(\tilde{n})f(x,1)=f( x-n, 1)\ .$

When $\tau=it$ is  purely imaginary, the representation $(\pi_\tau
,V_\tau)$ is pre-unitary and irreducible; the $G$-invariant scalar
product in $V_\tau$ is given by $\langle f,g \rangle_{V_\tau}=
\frac{1}{\pi}\int_{\br} f\bar g dx$. These representations are
called {\it the principal series} representations.

 When $\tau\in (-1,1)$, the representation $(\pi_\tau ,V_\tau)$ is called
 a representation of {\it the complementary series}. These representations
 are also pre-unitary and irreducible, but the formula for the scalar product is
 more complicated (see \cite{G5}).

All these representations have $K$-invariant vectors. We fix a
$K$-invariant unit vector $e_{\tau} \in V_\tau$ to be a function
which is constant on the unit circle $S^1$ in $\br^2$ in the plane
realization. Note that in the line model a $K$-fixed unit vector
is given by $e_\tau(x)=c(1+x^2)^{\frac{\tau-1}{2}}$ with
$|c|^2=\pi\inv$ for $\tau\in i\br$.

Another realization, which we call circle or spherical model, is
obtained by restricting functions in $V_\tau$ to the unit circle
$S^1\subset \br^2\sm 0$. In the circle model we have the
isomorphism $V_\tau\simeq C^\8_{even}(S^1)$ and for $\tau\in
i\br$, the scalar product is given by $\langle
f,g\rangle=\frac{1}{2\pi}\int_{S^1}f\bar gd\theta$ while the
action of $K$ is induced by the rotation of $S^1$.

Representations of the principal and the complimentary series
exhaust all nontrivial irreducible pre-unitary representations of
$G$ of class one.

\subsection{Automorphic representations} Every automorphic form
$\phi$ generates (under the right
translations by elements in $G$) an automorphic representation of
the group $G$ (see \cite{GGPS}); this means that, starting from
$\phi$, we produce a smooth irreducible unitarizable
representation of the group $G$ in a space $V$ and its realization
$\nu : V \to C^{\8}(X)$ in the space of smooth functions on the
automorphic space $X = \G \backslash G$. We will denote by
$V_\tau$ the isomorphism class of the representation arising in
this way from a Maass form $\phi=\phi_\tau$ with the eigenvalue
$\mu=\frac{1-\tau^2}{4}$.

Suppose we are given a class one  unitary representation $(\tau,
V_\tau)$ and an automorphic realization of it $\nu: V_{\tau} \to
C^\8(X)$; we assume $\nu$ to be an isometric embedding. Such $\nu$
gives rise to an eigenfunction of the Laplacian on the Riemann
surface $Y = X/K$ as before. Namely, if $e_{\tau} \in V_\tau$ is a
unit $K$-fixed vector, then the function $\phi = \nu(e_\tau)$ is a
$L^2$-normalized eigenfunction of the Laplacian on the space $Y =
X/K$ with the eigenvalue $\mu=\frac{1-\tau^2}{4}$. This explains
why $\tau$ is a natural parameter to describe Maass forms.

This construction gives a one-to-one correspondence between Maass
forms and class one automorphic representations (and more
generally between automorphic forms and automorphic
representations of $G$). We use this correspondence to translate
problems in automorphic forms into problems in representation
theory.

\section{Unipotent Fourier coefficients}\label{unipot-section}

\subsection{Whittaker functionals}\label{u-fourier} We start with
the classical interpretation of Fourier coefficients
$a_k(\phi_\tau)$ in terms of representation theory. Namely, we
consider Whittaker functionals on irreducible unitary
representations of $G$.

Let $\phi$ be a Maass form and $\nu:V\to C^\8(X)$ the
corresponding automorphic realization of the space of smooth
vectors of an irreducible unitary representation of $G$.

Let $N\subset G$ be the standard upper-triangular unipotent
subgroup. We denote by $\CN$ the $N$-invariant closed cycle (a
horocycle) $\G_\8\sm N\subset X$. The cycle $\CN$ is an orbit
$\CN=\bar e\cdot N\subset X$ of $N$, where  $\bar e$ is the image
of the identity element $e\in G$ under the natural projection
$G\to X$. In what follows we can choose any closed orbit of any
unipotent subgroup of $G$. We endow $\CN$ with the $N$-invariant
measure $dn$ of the {\it total mass one}, and fix an
identification $\G_\8\sm N\simeq\bz\sm\br$.

For $k\in\bz$, let $\psi_k:N\to\bc$ be the additive character
$\psi_k(t)=e^{2\pi ikt}$ of $N\simeq\br$ trivial on
$\G_\8\simeq\bz\subset\br\ $. We consider the functional
$l^a_k=l^{aut}_{\psi_k}:V\to\bc$ defined by the automorphic period
$$l^a_k(v)=\int_{\CN}\nu(v)(n)\bar\psi_{k}(n)dn\ \ {\rm for\ any\ } v\in
V\ .$$

The functional $l^a_k\in V^*$ is $(N,\psi_k)$-equivariant, i.e.,
$l^a_k(\pi(n)v)=\psi_k(n)l^a_k(v)$ for any $n\in N$ and $v\in V$.
It is well-known that for a {\it non-trivial} character $\psi_k$
the space of functionals in $V^*$ satisfying this property is
one-dimensional. The automorphic representation $(V,\nu)$ is
called cuspidal if $l^a_{\psi_0}\equiv 0$ (for any cuspidal
subgroup $\G_N$). We also have the standard Fourier expansion of
cuspidal automorphic functions along $\CN$:
\begin{eqnarray*}\nu(v)(x)=\sum_{k\not= 0} l^a_k(\pi(g)v),
\end{eqnarray*}
where $g$ corresponds to $x$ under the projection $p:G\mapsto
\G\sm G=X$ (i.e., $p(g)=x$).

 We now consider {\it model} Whittaker functionals.
In the line model of the representation $V=V_\tau\subset
C^\8(\br)$, we can construct a model Whittaker functional
$l_k^m=l_{\psi_k}^{mod}:V\to\bc$ by using the Fourier transform.
Namely, let $v\in V\subset C^\8(\br)$ be a vector (i.e., a smooth
function) of a compact support and $\xi\in\br$. We define the
model Whittaker functional by the integral
$$l^m_\xi(v)=\hat v(\xi)=\int_\br v(x)e^{-i\xi x}dx\ .$$
The functional $l^m_\xi$ clearly extends to the whole space $V$ by
continuity.

The uniqueness of Whittaker functionals implies that the model and
the automorphic functionals are proportional. Namely, for any
$k\in\bz\sm 0$, there exists a constant $a_k(\nu)\in\bc$ such that
\begin{eqnarray*}l^a_k=a_k(\nu)\cdot l^m_k\ .
\end{eqnarray*}

A simple computation shows that under our normalization $
|a_k(\nu)|=|a_k(\phi_\tau)|$. Namely, we have
$l^m_\xi(e_\tau)=\int(1+t^2)^{\frac{\tau-1}{2}}\exp(-i\xi
t)dt=\frac{\pi^\haf|\xi/2|^{-\tau/2}}{\G
(\frac{1-\tau}{2})}K_{-\tau/2}(\xi)$, where $K_t$ denotes the
$K$-Bessel function. Based on this we choose in \eqref{f-c} the
following normalization for the Whittaker functions
$$\CW_{\tau,k}(y)=l^{mod}_{\psi_k}\left(\pi_\tau\left(\begin{smallmatrix}y^\haf&\\
&y^{-\haf}\end{smallmatrix} \right)
e_\tau\right)=\frac{\pi^{\haf}}{\G
\left(\frac{1-\tau}{2}\right)}\cdot y^{-\haf}\
K_{-\tau/2}(2\pi|k|y)\ .$$ Under such normalization we have
$a_k(\phi)=a_k(\nu)$, and this is consistent with one of the
classical normalizations for Fourier coefficients of Maass forms
(see \cite{Iw}).

\subsection{Weighted sums of coefficients} We are interested in bounds
for Fourier coefficients $a_k(\nu)$. To this end we consider
bounds for weighted sums of the type
$$\sum_{k} |a_k(\nu)|^2\hat\al(k),$$ where $\hat\al$ is a
non-negative weight function. There is a simple geometric way to
construct these sums.

Let $\bar V$ be the representation which is complex conjugate to
$V$; it is also an automorphic representation with the realization
$\bar\nu:\bar V\to C^\8(X)$. We only consider the case of
representations of the principal series, i.e., we assume that $V =
V_{\tau}$, $\bar V= V_{-\tau}$ for some $ \tau \in i \mathbb{R}$;
\ the case of representations of the complementary series is
similar.

Consider the space $E=V\otimes \bar V$. We identify it with a
subspace of $C^\8(\br^2)$ using the line realization $V\subset
C^\8(\br)$. We have the corresponding automorphic realization
$\nu_E=\nu\otimes\bar\nu:E=V\otimes \bar V\to C^\8(X\times X)$.

Let $\Dl \CN\subset \Dl X\subset X\times X$ be the diagonal copy
of the cycle $\CN$.  We define the following automorphic
$N$-invariant functional $l_{\Dl \CN}:E\to\bc$ by
$$l_{\Dl \CN}(w)=\int_{\Dl \CN}\nu_E(w)(n,n)dn$$
for any $w\in E$.

We have the obvious Plancherel formula
\begin{equation}
\begin{split}\label{Fourier-l-N-w}
l_{\Dl \CN}(w)=\sum_k l^a_k\otimes \bar l^a_{-k}(w)=&\\
\sum_k|a_k(\nu)|^2 l^m_k\otimes \bar l^m_{-k}(w)=&
\sum_k|a_k(\nu)|^2 \hat w(k,-k)\
\end{split}\end{equation}
for any $w\in E\subset C^\8(\br^2)$ (here $\hat w$ denotes the
standard Fourier transform on $\br^2$).

Varying the vector $w\in E$, we obtain different weighted sums
$\sum_k|a_k(\nu)|^2\hat \al(k)$ with a weight function
$\hat\al(k)=\hat w(k,-k)$. The weight function might be easily
arranged to be non-negative as we will see below.

We now obtain another expression for the functional $l_{\Dl \CN}$
using spectral decomposition of $L^2(X)$ and trilinear invariant
functionals on irreducible representations of $G$. This will give
an instance of the Rankin-Selberg type formula discussed in
Introduction which in fact coincides with the classical formula of
Rankin and Selberg. We first discuss spectral decomposition of
$L^2(X)$ into irreducible unitary representations of $G$.

\subsection{Spectral decomposition and the Eisenstein series}
It is well-known that $L^2(X)=L^2_{cusp}(X)\oplus
L^2_{res}(X)\oplus L^2_{Eis}(X)$ decomposes into the sum of three
closed $G$-invariant subspaces  of cuspidal representations,
representations associated to residues of Eisenstein series and
the space generated by the unitary Eisenstein series (see
\cite{B}). The spaces $L^2_{cusp}(X)$ and $L^2_{res}(X)$ decompose
discretely into a direct sum of irreducible unitary
representations of $G$ and $L^2_{Eis}(X)$ is a direct integral of
irreducible unitary representations of the principal series. We
assume for simplicity that the residual spectrum is trivial, i.e.,
$L^2_{res}(X)=\bc$ is the trivial representation of $G$ (e.g.,
$\G$ is a congruence subgroup of $PSL_2(\bz)$).

We are interested in the spectral decomposition of the functional
$l_{\Dl \CN}$ defined as a period along the diagonal copy of a
horocycle inside of $X\times X$. Hence, the space $L^2_{cusp}(X)$
will not appear in the final formula as by the definition it
consists of functions satisfying $\int_\CN f(nx)dn=0$ for almost
all $x\in X$.

We will need the following basic facts from the theory of the
Eisenstein series (see \cite{Be}, \cite{B}, \cite{Ku}).

Let $B=AN$ be the Borel subgroup of $G$ (i.e., the subgroup of the
upper triangular matrices). We denote $\G_B=\G\cap B$,
$\G_N=\G_\8=\G\cap N$ and $\G_L=\G_B/\G_N$ which we assume, for
simplicity, is trivial. Let $Aff=N\sm G\simeq \{\br^2\sm 0\}/\{\pm
1\}$ be the basic affine space. The group $G$ acts from the right
on the space $Aff$ and preserves an invariant measure $\mu_{Aff}$.
The subgroup $B/N$ acts on $Aff$ on the left and acts on
$\mu_{Aff}$ by a character.

We denote $X_B=\G_BN\sm G$ and endow it with the measure
$\mu_{X_B}$ compatible with the measure $\mu_X$. We identify $X_B$
with $Aff$ (in general, one considers the space $\G_L\sm Aff$).

Let $\CA(X_B)$ be the space of smooth functions of moderate growth
on $X_B$. For a complex number $s\in \bc$, we denote by
$\CA^s(X_B)\simeq \CA^s(Aff)\simeq\CA^s_{even}(\br^2\sm 0)$ the
subspace of the even homogeneous functions of the homogeneous
degree $s-1$. The subspace $\CA^s(X_B)$ is $G$-invariant and for
$s$ purely imaginary is isomorphic to the space of smooth vectors
of a unitary class one representation of $G$.

In this setting one have the Eisenstein series operator
\begin{eqnarray*}
\E: \CA(X_B)\to C^\8(X)
\end{eqnarray*}
given by $\E(f)=\sum_{\g\in\G/\G_B}\g\circ f$ and the conjugate
constant term operator $ C: C^\8(X)\to \CA(X_B)\ $ given by
$C(\phi)=\int_{n\in N/\G_N}n\circ\phi\ dn$. The operator $\E$ is
only partially defined as the Eisenstein series not always
convergent.

The operators $\E$ and $C$ commute with the action of $G$. Hence
we also have the operator $\E(s)=\E|_{\CA^s(X_B)}:\CA^s(X_B)\to
C^\8(X)$ (defined via the analytic continuation for all $s\in
i\br$) and the fundamental relation $C(s)\circ \E(s)=Id+I(s)$
where $I(s):\CA^s(X_B)\to\CA^{-s}(X_B)$ is an intertwining
operator which is unitary for $s\in i\br$. It is customary to
write it in the form $I(s)=c(s)I_s$ where $I_s$ is a properly
normalized unitary intertwining operator satisfying $I_s\circ
I_{-s}=Id$ and $c(s)$ is a meromorphic function, satisfying the
functional equation $c(s)c(-s)=1$. The operator $I_s$ is
constructed explicitly in a model of the representation $V_s$. We
also have the functional equation $\E(s)=\E(-s)\circ I(s)$ for the
Eisenstein series.

The spectral decomposition of $L^2_{Eis}(X)$ then reads
$$L^2_{Eis}(X)=\int_{i\br^{+}}\E(s)(\CA^s(X_B))\ ds\ .$$ This
means, in particular, that for any $f\in C^\8(X)\cap L^2(X)$, the
Eisenstein component $f_{Eis}=pr_{Eis}(f)$ of $f$ in the space
$L^2_{Eis}(X)$ has the following representation $f_{Eis}=
\int_{i\br^+}\E(s)f_s\ ds$ for an appropriate smooth family of
functions $f_s\in \CA^s(X_B)$. We choose an orthonormal basis
$\{e_i(s)\}\subset \CA^s(X_B)$ and set $f_s=\sum_i \langle
f,\E(s)e_i(s)\rangle_{L^2(X)}e_i(s)$  for all $s\in i\br$. We have
then a more symmetrical spectral decomposition
$$f_{Eis}=
\haf\int_{i\br}\E(s)f_s\ ds\ ,$$ and the corresponding Plancherel
formula
$||f_{Eis}||_{L^2(X)}^2=\haf\int_{i\br}||f_s||_{\CA^s(X_B)}^2\
ds$.

\subsection{Trilinear invariant functionals}\label{triple} We construct the
spectral decomposition of $l_{\Dl \CN}$ with the help of trilinear
invariant functionals on irreducible unitary representations of
$G$. We review the construction below (for a more detailed
discussion see \cite{BR3}).

Let $\nu:V\to C^\8(X)$ be a cuspidal automorphic representation.
Let $E=V\otimes\bar V$ and $\nu_E$ be as above. Consider the space
$C^\8(X\times X)$. The diagonal $\Dl X \to X\times X$ gives rise
to the restriction morphism $r_\Dl : C^\8(X\times X)\to C^\8(X)$.
Let $\nu_W:W\to C^\8(X)$ be an irreducible automorphic
subrepresentation. We assume that for any $w\in W$ the function
$\nu_W(w)$ is a function of moderate growth on $X$. We define the
following $G$-invariant trilinear functional $l_{E\otimes
W}^{aut}=l_{\nu_E\otimes\nu_W}^{aut}$ on $E\otimes \bar W$ via
$$l_{E\otimes W}^{aut}(v\otimes v'\otimes u)=\langle r_\Dl(v\otimes v'),
u\rangle_{L^2(X)}$$ for any $v\otimes v'\in E$ and $u\in \bar W$.
The cuspidality of $V$ and the moderate growth condition on  $W$
ensure that $l_{E\otimes W}^{aut}$ is well-defined (i.e., the
integral over the non-compact space  $X$ is absolutely
convergent).

Next we use a general result from representation theory, claiming
that such a $G$-invariant trilinear functional is unique up to a
scalar (see \cite{O}, \cite{Pr} and the discussion in \cite{BR3}).
Namely, we have the following

\begin{thm}{ubi}Let $(\pi_j,V_j),$ where $j=1,2,3,$
be three irreducible smooth admissible
  representations of  $G$. Then
$\dim\Hom_G(V_1\otimes V_2\otimes V_3,\bc)\leq 1$.
\end{thm}

This implies that the automorphic functional $l^{aut}_{E\otimes
W}$ is proportional to an explicit \lq\lq model" functional
$l^{mod}_{E\otimes W}$ which we describe using explicit
realizations of representations $V$ and $W$ of the group $G$; it
is important that this model functional carries no arithmetic
information. The model functional is defined on any three
irreducible admissible representations of $\PGLR$ regardless
whether or not these are automorphic .

Thus we can write
\begin{eqnarray}\label{c-a}
l^{aut}_{E\otimes W} = a_{\nu_E\otimes \nu_W}\cdot
l^{mod}_{E\otimes W}
\end{eqnarray}
for some constant $a_{E\otimes W}=a_{\nu_E\otimes\nu_W}$ (somewhat
abusing notations as this coefficient depends on the realizations
$\nu_E$ and $\nu_W$ and not only on the isomorphism classes of $E$
and $W$).

It turns out that the proportionality coefficient $ a_{E\otimes
W}$ above carries important \lq\lq automorphic" information (e.g.,
essentially is equal to the Rankin-Selberg $L$-function) while the
second factor carries no arithmetic information and can be
evaluated using explicit realizations of representations $V$ and
$W$ (see Appendix in \cite{BR3} for an example of such a
computation). In what follows we only need the case of $W$ being
an irreducible unitary representation of the principal series
$V_s$, $s\in i\br$ (or the trivial representation).

Denote by $l^{mod}_s$ the {\it model} trilinear form $l^{mod}_s: V
\otimes \bar V \otimes V_s \to \bc$ which we describe explicitly
in Section \ref{mod-functional} below. Any $G$-invariant form $l:
V \otimes \bar V \otimes V_s \to \bc$ gives rise to a
$G$-intertwining morphism $T^l: V \otimes\bar V \to V_s^*$ which
extends to a $G$-morphism $ T^l: E \to \bar V_s$, where we
identify the complex conjugate space $\bar V_s$ with the smooth
part of the space $V_s^*$ ($\bar V_s\simeq V_{-s}$ for $s\in
i\br$). We apply this construction in order to describe the
projection of $E$ to the space orthogonal to cusp forms, namely to
$\bc\oplus L^2_{Eis}(X)$ (in general to $ L^2_{res}(X)\oplus
L^2_{Eis}(X)$).

We realize the irreducible principal series representation $V_s$
in the space of homogenous functions on the plane
$\CA^s(\{\br^2\sm 0\}/\{\pm 1\})\simeq\CA^s(Aff)\simeq
\CA^s(X_B)$. This is a model suitable for the theory of Eisenstein
series. For a chosen family of $G$-invariant functionals
$l_s^{mod}=l_{E\otimes V_s}:E\otimes V_{-s}\to\bc$ and the
corresponding family of morphisms $T_s=T^{l_s^{mod}}:E\to
V_{s}\simeq A^s(X_B)$, we have the proportionality coefficient
$a(s)=a_\nu(s)=a_{E\otimes V_{-s}}$ defined by $l^{aut}_s=a(s)\
l^{mod}_s$ as in \eqref{c-a} and the corresponding spectral
decomposition for any $w\in E$,
\begin{align}\label{E-Eis-proj-sym}
\ pr_{res\oplus Eis}(\nu_E(w))=\langle r_\Dl(\nu_E(w)),
1\rangle\cdot 1+ \haf\int_{i\br}a(s)\E(s)(T_s(w))\ ds\ .
\end{align}
Note that \eqref{E-Eis-proj-sym} is symmetrical under the change
$s\to -s$. This is achieved by  choosing the model trilinear
functionals $l^{mod}_{ s}:E\otimes V_{ s}\to \bc$ to satisfy the
relation $l^{mod}_s=l^{mod}_{-s}\circ I_s$ and the coefficients
$a(s)$ to satisfy $a(s)=c(s)a(-s)$ (this is equivalent to the
functional equation for the Rankin-Selberg $L$-function). We note
also that $\langle r_\Dl(\nu_E(w)), 1\rangle=Tr(w)$ for any $w\in
E$, viewed as an element in $V\otimes V^*$.

We use spectral decomposition \eqref{E-Eis-proj-sym} to obtain the
spectral decomposition of the functional $l_{\Dl\CN}$.

Namely, consider the functional $l_\CN:C^\8(X)\to\bc$ given by the
constant term along $\CN\subset X$ (i.e.,
$l_{\Dl\CN}(f)=C(f)|_{x=\bar e}$ for any $f\in C^\8(X)$). As
$l_\CN$ vanishes on $L^2_{cusp}(X)$, we only have to understand
its form on the space of the Eisenstein series (and on the space
of residues). The pair $(G,N)$ is not a Gelfand pair (the space of
$N$-{\it invariant}  functionals is two dimensional) and we can
not use the argument we used for the Whittaker functionals.
However, the theory of the Eisenstein series provides the
necessary remedy. Namely, consider the representation of the
(generalized) principal series $\CA^s$ realized in the space of
even homogenous functions on $X_B\simeq \br^2\setminus 0$. The
space of $N$-invariant functionals on $\CA^s$ is generated by the
functionals $\dl_s$ and $\dl_{-s}$, where $\dl_s(v)=v(0,1)$ and
$\dl_{-s}(v)=I_s(v)(0,1)$ (in fact, the functional $\dl_{-s}$ is
given (up to a normalization constant) by the integral over the
line $\{(1,x)|\ x\in\br\}\subset\br^2$). The basic theory of the
constant term of the Eisenstein series then implies that
\begin{eqnarray*}\label{const-term}
C(\E(s)(v))|_{x=\bar e}=\dl_s(v)+c(s)\dl_{-s}(v)\ .
\end{eqnarray*}

Applying this to \eqref{E-Eis-proj-sym}, we obtain the following
spectral decomposition
\begin{eqnarray*}\label{l-N-RS}
l_{\Dl\CN}(\nu_E(w))=l_\CN(pr_{res\oplus
Eis}(\nu_E(w)))=\frac{vol(\CN)}{vol(X)^{-\haf}}\cdot Tr(w)+
\int_{i\br}a(s)\dl_s(T_s(w))\ ds\ .
\end{eqnarray*}
Here we have used the functional equation for the constant term of
the Eisenstein series
\begin{eqnarray*}
a(s)c(s)\cdot\dl_{-s}(T_{s}(w))= a(-s)\cdot\dl_{-s}(T_{-s}(w))\
\end{eqnarray*} and the assumption that the residual spectrum is
trivial. Taking into consideration the Plancherel formula
\eqref{Fourier-l-N-w} and the normalization of measures
$vol(X)=vol(\CN)=1$, we arrive at the identity
\begin{eqnarray}\label{RS2}
\sum_k|a_k(\nu)|^2 \hat w(k,-k)= Tr(w)+
\int_{i\br}a(s)\dl_s(T_s(w))\ ds\ .
\end{eqnarray}

This is our form of the Rankin-Selberg formula. To give it a more
familiar form similar to \eqref{RS-Fourier}, we will make
\eqref{RS2} more explicit by describing $T_s$ and $\dl_s$ in the
line model of $V_s$. We do this by choosing an explicit kernel for
the model invariant trilinear functional $l^{mod}_s$.

\subsection{Model trilinear functionals}\label{mod-functional} We recall
the construction of model trilinear functionals presented in
\cite{BR3}. There it was shown that in the line model of
representations $V\simeq V_\tau$, $\bar V\simeq V_{-\tau}$ and
$V_{-s}$ the kernel
\begin{equation}\begin{split}
\label{kernel} K_{\tau,-\tau,-s}(x,y,z)=&\\ |x-y|^{(-s-1)/2}&
|xz-1|^{(-2\tau+s-1)/2}|yz-1|^{(2\tau+s-1)/2}\
\end{split}\end{equation} defines a nonzero
trilinear $G$-invariant functional $l^{mod}_s$ on $V_\tau\otimes
V_{-\tau}\otimes V_{-s}$. This gives rise to the map
$T_s:E_\tau\simeq V_\tau\otimes V_{-\tau} \to V_s$ given by the
same kernel. Here the variable $x$ corresponds to the
representation $V_\tau$ and $y,\ z$ to $V_{-\tau}$ and $ V_{-s}$,
respectively. Note a certain asymmetry between $V_\tau,\
V_{-\tau}$ and $V_{-s}$. This is because we choose for the first
two representations the line model associated with the upper
triangular subgroup and for the last representation the model
associated with the lower triangular subgroup.

In the line model $V_s\subset C^\8(\br)$ the $N$-invariant
functional $\dl_s$ is given by the evaluation at the point $z=0$:
$\dl_s(f)=f(0)$. Hence from \eqref{kernel} it follows that the
composition $\dl_s(T_s)$ is given by the Mellin transform for any
$w\in E\subset C^\8(\br\times \br)$,
\begin{eqnarray*}\label{melin-E}
\dl_s(T_s(w))=\int_{\br^2} w(x,y)|x-y|^{(-s-1)/2}dxdy\ .
\end{eqnarray*}

Plugging this into \eqref{RS2}, we arrive at the ``classical"
Rankin-Selberg formula (we assume as before that the residual
spectrum is trivial)
\begin{eqnarray}\label{RS4}
\sum_k|a_k(\nu)|^2 \hat w(k,-k)=Tr(w)+ \int_{i\br}a(s)w^\flat(s)\
ds\ ,
\end{eqnarray}
where we denoted
\begin{eqnarray}\label{flat-int} w^\flat(s)=\int
w(x,y)|x-y|^{(-s-1)/2}dxdy\ . \end{eqnarray} This coincides with
the Mellin transform $M(\al)(s)$ for
$\al(t)=\int_{x-y=t}w(x,y)dl$. The transform $^\flat$ is
 defined for any smooth rapidly decreasing function $w$, at
least for all ${\rm Re}(s)< 1$. In fact, it could be defined for
all $\lm\in \bc$, by means of analytic continuation, but we will
not need this. We will consider only the case $s\in i\br$ as we
assume that the residual spectrum is trivial. In general, the
residual spectrum could be treated similarly. We note also that
$Tr(w)=\int w(x,x)dx=\al(0)$.

We can now rewrite  the Rankin-Selberg formula in a more familiar
form if we substitute the vector $w(x,y)$ by
$\al(t)=\int_{x-y=t}w(x,y)dl$. We have
\begin{eqnarray}\label{RS3}
\sum_k|a_k(\nu)|^2 \hat\al(k) =\al(0)+ \int_{i\br}a(s)\cdot
M(\al)(s)\ ds\ ,
\end{eqnarray}
where $\hat\al(\xi)=\hat w(\xi,-\xi)$. This coincides with the
classical Rankin-Selberg formula.

\subsubsection{Remarks}\label{RS-Fourier-Mellin-form} 1.
Taking into account that $M(\al)(s)=\g(s)M(\hat \al)(1-s)$, where
$\g(s)=
\frac{\pi^{-\frac{s}{2}}\G(\frac{s}{2})}{\pi^{-\frac{1-s}{2}}\G(\frac{1-s}{2})}$
(note that $|\g(s)|=1$ for $s\in i\br$), we see that
\begin{eqnarray}\label{RS5}
\sum_k|a_k(\nu)|^2 \hat\al(k) =\al(0)+
\int_{i\br}a(s)\g(s)M(\hat\al)(s)\ ds\ .
\end{eqnarray}
This seems to have an advantage of being an identity for one
function $\hat\al$ and not for two functions $\hat\al$ and
$M(\al)$. In practice we find it easier to work with one master
function $\al$ and use the identity \eqref{RS3}.

2. We would like to point out one essential difference between the
classical Rankin-Selberg formula \eqref{RS-Fourier} obtained via
the unfolding and formula \eqref{RS3} that we prove. The unfolding
method provides an explicit relation between a choice of a model
Whittaker functional on a cuspidal representation and the
coefficient of proportionality $D(s,\phi,\bar\phi)$ (i.e.,
essentially the Rankin-Selberg $L$-function). In the argument we
presented, the coefficient of proportionality $a(s)$ in addition
depends on the choice of the auxiliary model trilinear functional.
One can use the Whittaker functional in order to define the model
trilinear functional and hence eliminate this extra indeterminacy.
We hope to return to this subject elsewhere.

\subsection{Proof of Theorem \ref{thm}}\label{proof-F-u}
We are interested in getting a bound for the coefficients
$a_n(\phi)$. The idea  is to find a test vector $w\in V\otimes
\bar V$, i.e., a function $w\in C^\8(\br\times\br)$, such that
when substituted in the Rankin-Selberg formula \eqref{RS4} will
produce a weight $\hat w$ which is not too small for a given $n$
($|n|\to\8$). Given such a vector we have to estimate its spectral
density, i.e., the transform $w^\flat$. One might be tempted to
take $w$ such that $\hat w$ is essentially a delta function (i.e.,
the weight $\hat w$ picks up just a few coefficient $a_n(\phi)$ in
\eqref{RS4}). However, for such a vector we have no means to
estimate the right hand side of the Rankin-Selberg formula
\eqref{RS4} because the weight function $w^\flat$ is spread over
an interval of the spectrum which is too long to use the maximum
modulus bound (still, conjecturally even for such functions  the
contribution on the right hand side of the Rankin-Selberg formula
is small thanks to cancellations). The solution to this problem is
well-known in harmonic analysis. One takes a function which
produces a weighted sum of the coefficients $|a_k(\phi)|^2$ in a
certain range depending on $n$ and such that its transform
$w^\flat$ spreads over a shorter interval. For a certain kind of
such test vectors $w$ (namely, those for which the support of
$\hat w$ is not too small), we give essentially a {\it sharp}
bound for the value of $l_{\Dl\mathcal{N}}(w)$.

We now explain how to choose the required test vectors. Let $\chi$
be a smooth function on $\br$ with the support
$supp(\chi)\subset[-\haf,\haf]$ and such that the Fourier
transform satisfies $|\hat\chi(\xi)|\geq 1$ for $|\xi|\leq 1$. We
consider the convolution $\psi=\chi*\chi'$, where
$\chi'(x)=\bar\chi(-x)$. We have $supp(\psi)\subset[-1,1]$,
$\hat\psi(\xi)\geq 0$ for all $\xi$ and $\hat\psi(\xi)\geq 1$ for
$|\xi|\leq 1$.

Let $N\geq T\geq 1$  be two real numbers. We consider the
following test vector
\begin{eqnarray*}
w_{N,T}(x,y)=T\cdot e^{-iN(x-y)}\cdot\psi(T(x-y))\cdot \psi(x+y)\
.
\end{eqnarray*}
We have the following technical lemma describing properties of
$w_{N,T}^\flat$ (where the transform $^\flat$ was defined in
formula \eqref{flat-int} and is essentially the Mellin transform
in $|x-y|$).

\begin{lem}{lem1}For $w_{N,T}$ as above, the following bounds
hold:
\begin{enumerate}
    \item $|\int w_{N,T}(t,t)dt|\leq c T$,
    \item $\hat w_{N,T}(\xi,-\xi)\geq 0$ for all $\xi$,
    \item $\hat w_{N,T}(\xi,-\xi)\geq 1$ for all $\xi$ such that $|\xi-N|\leq T$,
    \item $|w_{N,T}^\flat(s)|\leq c
    T|N|^{-\haf}$ for $|s|\leq N/T$,
    \item $|w_{N,T}^\flat(s)|\leq c
    T(1+|s|)^{-3}$ for $|s|\geq N/T$,
\end{enumerate} for some fixed constant $c>0$ which is independent of $N$ and $T$.
\end{lem}

Bounds $(1)$-$(3)$ are immediate. Bounds $(4)$ and $(5)$ are
standard, once we apply the stationary phase method or the van der
Corput lemma to the integral of the form
$w_{N,T}^\flat(s)=\hat\psi(0)\cdot T^{\haf+s/2}\cdot\int \psi(t)\
e^{-i\frac{N}{T}t}|t|^{-\haf-s/2}dt$ (see Section
\ref{proof-lem1}).

We return to the proof of Theorem \ref{thm-2}. We will use the
following mean value (or convexity) bound
\begin{eqnarray}\label{bound-1}
\int_0^A|a(it)|^2dt\leq C_\tau A^2\ln A\ ,
\end{eqnarray} proved in
\cite{BR1} for any $A\geq 1$. The constant $C_\tau$ satisfies the
bound $C_\tau\leq C_\G(1+|\tau|)$ with a constant $C_\G$ depending
on $\G$ only.

We substitute the vector $w_{N,T}$ into the Rankin-Selberg formula
\eqref{RS4} (and use $Tr (w)=\int w(t,t)dt$). Taking into account
the convexity bound \eqref{bound-1}, bounds in the lemma and the
Cauchy-Schwartz inequality, we obtain
\begin{align*}
\nonumber
 & \sum_{|k-N|\leq T}|a_k(\nu)|^2\leq\sum_k|a_k(\nu)|^2 \hat
  w_{N,T}(k)=\\ &
  =\ \int w_{N,T}(t,t)\ dt
  +\int_{i\br}a(s) w_{N,T}^\flat(s)\ d|s|\leq\\\nonumber
&\leq cT+\int_{|s|\leq N/T}cT|N|^{-\haf}|a(s)|\ d|s|+
\int_{|s|\geq N/T}cT(1+|s|)^{-3}|a(s)|\ d|s| \leq\\\nonumber &\leq
cT+cT|N|^{-\haf}\left(\int_{|s|\leq
N/T}|a(s)|^2d|s|\cdot\int_{|s|\leq N/T}1\ d|s|\right)^\haf\ +
\\\nonumber
 & \qquad\qquad\qquad +\ cT\int_{|s|\geq N/T}(1+|s|)^{-3}(1+|a(s)|^2)\ d|s|\leq \\\nonumber
 &\leq
  cT+CT|N|^{-\haf}\left(\frac{N}{T}\right)^{3/2+\eps}+DT
  =c'T+CT^{-\haf-\eps}|N|^{1+\eps}\ ,
\end{align*} for any $\eps>0$  and some constants $c',\ C,\ D>0$.

Setting $T=N^{2/3}$, we obtain $ \sum\limits_{|k-N|\leq
N^{2/3}}|a_k(\nu)|^2\leq A_\eps N^{2/3+\eps}$ for any $\eps>0$.
\qed

\subsection{Remarks}

1. It is more customary to use the formula \eqref{RS5}. We find
the geometric formula \eqref{RS4} more transparent. Following the
argument of Good \cite{Go1}, one usually argues as follows. For
$R\geq 1$ and $Z\geq 1$, choose a test function
$\al_{Z,R}(t)=\al_Z(t/R)$, where $\al_Z$ is smooth, supported in
$(1-2/Z,1+2/Z)$ and $\al|_{(1-1/Z,1+1/Z)}\equiv 1$. This means
that the sum $\sum_k|a_k(\nu)|^2 \al_{Z,R}(k)$ is essentially over
$k$ in the interval of size $R/Z$ centered at $R$. The Mellin
transform $M(\al_Z)(s)=\int_{\br^+}\al_Z(t)|t|^{s}d^xt$ of $\al_Z$
satisfies the simple bound
$$|M(\al_Z)(s)|\leq cZ\inv$$ for any $|s|$, and the bound
$$|M(\al_Z)(s)|\leq c|s|\inv\left(\frac{Z}{|s|}\right)^{m}$$ for
any $m>0$ and $|s|\geq 1$. This follows easily  from integration
by parts (we are only interested in $s\in i\br$). In particular,
we have $|M(\al_Z)(s)|\leq cZ^{\haf+\eps}|s|^{-3/2-\eps}$ for
$|s|\geq Z$. Using the average bound $ \int_0^A|a(it)|^2dt\leq
C_\tau A^2\ln A $ and  the Cauchy-Schwartz inequality, one obtains
\begin{eqnarray*}\label{bound-2-0}
\left|\int_{i\br} a(s)\g(s)M(\al_{Z,R})(s)ds\right|\leq C_\eps
R^{\haf+\eps} Z^{\haf+\eps}
\end{eqnarray*} for any
$\eps>0$. We arrive at $\sum_k|a_k(\nu)|^2 \al_{Z,R}(k) \leq  R/Z+
C_\eps R^{\haf+\eps} Z^{\haf+\eps}$ and setting $Z=R^{1/3}$, we
obtain the bound claimed.

2. One might conjecture that for any $A\geq 1$, the bound
\begin{eqnarray}\label{Lindlef}
\int_A^{2A}|a(it)|^2dt\ll_{\nu,\eps} A^{1+\eps} \end{eqnarray}
holds for any $\eps>0$ (e.g., the Lindel\"{o}ff conjecture on the
average for the Rankin-Selberg $L$-function). It is easy to see
that such a bound would lead to the bound
$|a_n(\nu)|\ll_{\nu,\eps} |n|^{\qtr+\eps}$. We note that this
bound is a natural barrier which for the Rankin-Selberg method
would be hard to overcome. Nevertheless, it was suggested by
Sarnak \cite{Sa3} that for a general lattice $\G\subset\PGLR$ the
Ramanujan-Petersson conjecture $|a_n(\phi_\tau)|\ll |n|^\eps$
might hold.

\subsection{Proof of Lemma \ref{lem1}}\label{proof-lem1}
We prove the following statement from which  Lemma \ref{lem1}
immediately follows.

\begin{lem}{lem3} Let $\psi$ be a fixed smooth function
with a compact support in $[-1,1]$. For $s\in i\br$ and
$\xi\in\br$, denote $\psi^\flat(\xi, s)=\int_\br \psi(t)e^{-i\xi
t}|t|^{-\haf-s}dt$. There exists a constant $c>0$ such that
\begin{enumerate}
    \item $|\psi^\flat(\xi, s)|\leq c
    (1+|\xi|)^{-\haf}$ for $|s|\leq 2|\xi|$,
    \item $|\psi^\flat(\xi, s)|\leq c
    (1+|s|)^{-3}$ for $|s|\geq 2|\xi|$.
\end{enumerate}
\end{lem}

In fact, both claims in the lemma are simple consequences of the
van der Corput lemma (see \cite{St}, p.~332). One also can use the
following Fourier transform argument.

To prove $(1)$, we use the fact that the Fourier transform of
$|t|^{-\haf-s}$ is equal to $\g(-\haf-s)|\xi|^{-\haf+s}$, where
$\g(s)$ is defined in Remark \ref{RS-Fourier-Mellin-form}, and
$|\g(-\haf-it)|=1$. The Fourier transform of $\psi$ satisfies
$|\hat\psi(\xi)|\ll (1+|\xi|)^{-M}$ for any $M>0$. Hence, the
Fourier transform of $\psi(t)|t|^{-\haf-s}$, i.e., the convolution
$\hat\psi(\xi)*|\xi|^{-\haf+s}$,  is bounded by
$c(1+|\xi|)^{-\haf}$ for some $c$ and all $s\in i\br$. This proves
$(1)$.

To prove $(2)$, it is enough to notice that under the condition
$|s|\geq|\xi|$ the phase in the oscillating integral defining
$\psi^\flat(\xi, s)$ has no stationary points. The resulting bound
easily follows from the integration by parts (see Appendix
\ref{apen-A} for the similar computation).\qed

\section{Spherical Fourier coefficients}\label{sperical-section}

When dealing with the spherical Fourier coefficients we assume,
for simplicity, that the lattice $\G$ is co-compact. A general
finite co-volume lattice could be treated analogously without any
significant changes (see Remark \ref{rem-sph}).
\subsection{Geodesic circles}\label{K-geometry}We start with
the geometric origin of the spherical Fourier coefficients. We fix
a maximal compact subgroup $K=PO(2)\subset G$ and the
identification $G/K\to\uH$, $g\mapsto g\cdot i$. Let $y\in Y$ be a
point and $p:\uH\to \G\sm \uH\simeq Y$ the corresponding
projection compatible with the distance function ${\rm
d}(\cdot,\cdot)$ on $Y$ and on $\uH$. Let $R_y>0$ be the
injectivity radius of $Y$ at $y$. For any $r< R_y$ we define the
geodesic circle of radius $r$ centered at $y$ to be the set
$\s(r,y)=\{y'\in Y|{\rm d}(y',y)=r\}$. Since the map $p$ is a
local isometry, we have that $p(\s_\uH(r,z))=\s(r,y)$ for any
$z\in\uH$ such that $p(z)=y$, where $\s_\uH(r,z)$ is the
corresponding geodesic circle in $\uH$ (all geodesic circles in
$\uH$ are the Euclidian circles, though with a center different
from $z$). We associate to any such circle on $Y$ an orbit of a
compact subgroup on $X$. Namely, let $K_0=PSO(2)\subset K$ be the
connected component of $K$. Any geodesic circle on $\uH$ is of the
form $\s_\uH(r,z)=hK_0g\cdot i$ with $h,\ g\in G$ such that
$h\cdot i=z$ and $hg\cdot i\in \s_\uH(r,z)$ (i.e. an
$h$-translation of a standard geodesic circle centered at
$i\in\uH$ and passing through $g\cdot i\in\uH$). Note, that the
radius of the circle is given by the distance ${\rm d}(i,g\cdot
i)$ and hence $g\not\in K_0$ for a nontrivial circle. Given the
geodesic circle $\s(r,y)\subset Y$, we choose a circle
$\s_\uH(r,z)\subset \uH$ projecting onto $\s(r,y)$ and the
corresponding elements $g,\ h\in G$. We denote by $K_\s=g\inv
K_0g$ the corresponding compact subgroup and consider its orbit
$\ck_\s=hg\cdot K_\s\subset X$. Clearly we have $p(\ck_\s)=\s$. We
endow the orbit $\ck_\s$ with the unique $K_\s$-invariant measure
$d\mu_{\ck_\s}$ of the {\it total mass one} (from the geometric
point of view a more natural measure would be the length of $\s$).

We note that in what follows the restriction $r<R_y$ is not
essential. From now on we assume that $\ck\subset X$ is a
connected orbit of a {\it connected} compact subgroup $K'\subset
G$ (i.e., $K'$ is conjugated to $PSO(2)$). The restriction $r<R_y$
simply means that the projection $p(\ck)\subset Y$ is a smooth
non-self intersecting curve on $Y$. We also remark that it is
well-known that polar geodesic coordinates $(r,\theta)$ centered
at a point $z_0\in\uH=G/K$ could be obtained from the Cartan
$KAK$-decomposition of $G$ (see \cite{He}). This allows one to
give a purely geometric construction of the functions $P_{n,\tau}$
from the Introduction (see Section \ref{P-n}).

\subsection{ $K'$-equivariant
functionals}\label{K-fourier}  We fix a point $\dot{o}\in\ck$. Let
$\chi:K'\to S^1$ be a character. To such a character we associate
a function $\chi_.(\dot{o}k')=\chi(k')$, $k'\in K'$ on the orbit
$\ck$ and the corresponding functional on $C^\8(X)$ given by
\begin{align*}\label{O-funct-K}
d^{aut}_{\chi,\ck}(f)=\int_\ck f(k)\bar{\chi}_.(k)d\mu_\ck
\end{align*}
for any $f\in C^\8(X)$.

The functional $d^{aut}_{\chi,\ck}$ is $\chi$-equivariant:
$d^{aut}_{\chi,\ck}(R(k')f)=\chi(k')\cdot d^{aut}_{\chi,\ck}(f)$
for any $k'\in K'$, where $R$ is the right action of $G$ on the
space of functions on $X$. For a given orbit $\ck$ and a choice of
a generator $\chi_1$ of the cyclic group $\hat{K'}\simeq\bz$ of
characters of the compact group $K'$, we will use the shorthand
notation $d^{aut}_{n}=d^{aut}_{\chi_n,\ck}$, where
$\chi_n=\chi_1^n$. The functions $\chi_{n}.\ $, $n\in\bz$ form an
orthonormal basis for the space $L^2(\ck,d\mu_\ck)$ (since we
normalized the measure by $\mu_\ck(\ck)=1$).

 Let $\nu:V\to C^\8(X)$ be an irreducible automorphic
representation. When it does not lead to confusion, we denote by
the same letter the functional
$d^{aut}_{\chi,\ck}=d^{aut}_{\chi,\ck,\nu}$ on the space $V$
induced by the functional $d^{aut}_{\chi,\ck}$ defined above on
the space $C^\8(X)$. Hence we obtain an element in the period
space $\mathcal{P}_{K'}(V,\chi)=\Hom_{K'}(V,\chi)$. We next use
the well-known fact that this space is at most one-dimensional.

Let $V\simeq V_\tau$ be a representation of the generalized
principal series. We have then $\dim\Hom_{K'}(V_\tau,\chi)\leq 1$
for any character $\chi$ of $K'$ (i.e., the space of $K'$-types is
at most one dimensional for a maximal connected compact subgroup
of $G$). In fact, $\dim\Hom_{K'}(V_\tau,\chi_n)=1$ if and only if
$n$ is even.

To construct a model $\chi$-equivariant functional on $V_\tau$, we
consider the circle model $V_\tau\simeq C^\8_{even}(S^1)$ in the
space of even functions on $S^1$ and the standard vectors
(exponents) $e_n=\exp( in\theta)\in C^\8(S^1)$ which form a basis
of $K$-types for the {\it standard} maximal compact subgroup
$K=PO(2)$. For any $n$ such that
$\dim\Hom_{K_0}(V_\tau,\chi_n)=1$, the vector
$e'_n=\pi_\tau(g\inv)e_n$ defines a non-zero
$(\chi_n,K')$-equivariant functional on $V_\tau$ by the formula
\begin{eqnarray*}
d^{mod}_{n}(v)=d^{mod}_{\chi_n,\tau}(v)=\langle v,e'_n\rangle\ .
\end{eqnarray*} We call such a functional the {\it model}
$\chi_n$-equivariant functional on $V\simeq V_\tau$.

The uniqueness principle implies that there exists a constant
$b_n(\nu)=b_{\chi_n,\ck}(\nu)$ such that for any $v\in V$
\begin{eqnarray*}\label{b-def}
d^{aut}_{n}(v)=b_n(\nu)\cdot d^{mod}_{n}(v)\ .
\end{eqnarray*}

\subsubsection{Functions $P_{n,\tau}$}\label{P-n}
We want to compare the coefficients $b_n({\nu})$ to the
coefficients $b_n(\phi_\tau)$ we introduced in
\eqref{spherical-expansion}. In particular we describe the
functions $P_{n,\tau}$ and their normalization. Let $h,\ g\in G$
and $\ck=hgK'\subset \G\sm G=X$ be the orbit of the connected
compact group $K'=g\inv K_0g$ as above. Let $\nu:V_{\tau}\to
C^\8(X)$ be an automorphic realization and
$\phi_{\tau}=\nu(e_0)\in C^\8(X)$ the $K$-invariant vector which
corresponds to a $K$-invariant vector $e_0\in V_{\tau}$ of norm
one, i.e., $\phi_{\tau}$ is a Maass form. We define the function
$P_{n,\tau}$ through the following matrix coefficient:
$P_{n,\tau}(r)e^{in\theta}=\langle e_0,\pi_\tau(g\inv
k\inv)e_n\rangle_{V_\tau}$, where $(r,\theta)=z=hkg\cdot i\in\uH$
for $k\in K_0$. It is well-known that the matrix coefficient is an
eigenfunction of the Casimir operator and hence
$P_{n,\tau}(r)e^{in\theta}$ is an eigenfunction of $\Dl$ on $\uH$.
In fact, the functions $P_{n,\tau}$ are equal to the Legendre
functions for the special choice of parameters (compare to
\cite{Iw}). Under such normalization of functions $P_{n,\tau}$, we
have
\begin{eqnarray*}
|b_n({\nu})|=|b_n(\phi_\tau)|\ .
\end{eqnarray*}

Let $\bar V$ be the complex conjugate representation; it  is also
an automorphic representation with the realization $\bar\nu:\bar
V\to C^\8(X)$. We only consider the case of representations of the
principal series, i.e., we assume that $V = V_{\tau}$, $\bar V=
V_{-\tau}$ for some $ \tau \in i \mathbb{R}$; \ the case of
representations of the complementary series can be treated
similarly. Let $\{e_n\}_{n\in 2\bz}$ be a $K$-type orthonormal
basis in $V$. We denote by $\{\bar e_{n}\}$ the complex conjugate
basis in $\bar V$ and denote by $\bar d^{aut/mod}_{n}$ the
corresponding automorphic/model functionals on the conjugate space
$\bar V\simeq V_{-\tau}$.

We introduce another notation for a $K'$-{\it invariant}
functional on an irreducible automorphic representation
$\nu_i:V_{\lm_i}\to C^\8(X)$ of class one.

Let $\chi_0:K'\to 1\in S^1\subset\bc$ be the trivial character of
$K'$. We have as above
\begin{eqnarray*}
d^{aut}_{\chi_0,\ck,\nu_i}(v)=\int_\ck
\nu_i(v)(k)\bar{\chi_0}_.(k)d\mu_\ck=b_0(\nu_i)\cdot \langle
v,e'_0\rangle_{V_{\lm_i}}\ ,
\end{eqnarray*}
for any $v\in V_{\lm_i}$.

We denote by $d_\lm(v)=\langle v,e'_0\rangle_{V_\lm}$ the
corresponding model functional and by
$$\beta(\lm_i)=b_0({\nu_i})$$ the proportionality coefficient
(somewhat abusing notations, since the coefficient depends on the
automorphic realization $\nu_i$ and not only on the isomorphism
class  $V_{\lm_i}$).

We want to compare coefficients $\beta(\lm_i)$ with more familiar
quantities. Let $\ck=x_0K'\subset X$ be an orbit of the compact
group $K'$. Let $\nu_i:V_{\lm_i}\to C^\8(X)$ be an automorphic
realization and $\phi'_{\lm_i}=\nu_i(e'_0)$ the $K'$-invariant
vector which corresponds to a $K'$-invariant vector $e'_0\in
V_{\lm_i}$ of norm one. From the definition of $b_0({\nu_i})$ it
follows that
\begin{eqnarray}\label{b-0-nomalize}
\beta(\lm_i)=\phi'_{\lm_i}(x_0)\ .
\end{eqnarray}
Hence, under the normalization we choose, the coefficients
$\beta(\lm_i)$ coincide with the value at a point $x_0$ for Maass
form $\phi'_{\lm_i}$ on the Riemann surface $Y'=\G\sm G/g\inv Kg$.

Finally, we note that on the discrete series representations any
$K'$-invariant functional is identically zero. This greatly
simplifies the technicalities in what follows.

\subsection{$\Dl\CK$-restriction}
Let $\Dl \ck\subset \Dl X\subset X\times X$ be the diagonal copy
of the cycle $\ck$.  We define the $\Dl K'$-invariant automorphic
functional $d_{\Dl \ck}:E=V\otimes\bar V\to\bc$ by
$$d_{\Dl \ck}(w)=\int_{\Dl \ck}\nu_E(w)(k,k)d\mu_\ck$$
for any $w\in E$.

Arguing as in Section \ref{u-fourier}, we have the following
Plancherel formula on $\ck$
\begin{equation}\begin{split}
\label{Fourier-d-K-w} d_{\Dl \ck}(w)=\sum_n d^{aut}_n\otimes \bar
d^{aut}_{-n}(w)&=\\ =\ \sum_n|b_n(\nu)|^2 & d^{mod}_n\otimes \bar
d^{mod}_{-n}(w)= \sum_n|b_n(\nu)|^2 \hat w(n,-n)\ ,
\end{split}
\end{equation}
where $\hat w(n,-n)=\langle w, e_n\otimes\bar e_{-n}\rangle_E$.
Hence, choosing various vectors $w$  we can obtain a variety of
weighted sums $\sum_n|b_n(\nu)|^2\hat\al(n)$.

We now obtain another expression for the functional $d_{\Dl \ck}$
using the spectral decomposition of $L^2(X)$ and trilinear
invariant functionals as in Section \ref{triple}.

\subsection{Anisotropic Rankin-Selberg formula. Proof of Theorem \ref{thm-2}}
 We assume that the space $X$ is compact. This
implies the discrete sum decomposition $L^2(X)=\left(\oplus_i
(L_i,\nu_i)\right)\oplus\left(\oplus_\kappa
(L_\kappa,\nu_\kappa)\right)$   into irreducible unitary
representations of $G$, where $\nu_i:L_i\simeq L_{\lm_i}\to
L^2(X)$ are unitary representations of class one (i.e., those
which correspond to Maass forms on $Y$) and $L_\kappa$ are
representations of discrete series (i.e., those which correspond
to holomorphic forms on $Y$). We fix such a decomposition and
denote by $V_i\subset L_i$ the corresponding spaces of smooth
vectors, and by $\nu_i^*:C^\8(X)\to V_i$ the adjoint map.

We fix $\nu:V\to C^\8(X)$ an irreducible automorphic
representation and denote by $\nu_E=\nu\otimes\bar\nu: E=V\otimes
\bar V\to C^\8(X\times X)$ the corresponding realization of $E$.

We use the notations from Section \ref{triple}. Let
$r_\Dl:C^\8(X\times X)\to C^\8(X)$ be the map induced by the
imbedding $\Dl: X\hookrightarrow X\times X$. Let
$\nu_i:V_{\lm_i}\to C^\8(X)$ be an irreducible automorphic
representation. Composing $r_\Dl$ with the adjoint map
$\nu_i^*:C^\8(X)\to V_{\lm_i}$, we obtain the trilinear $\Dl
G$-invariant map $T^{aut}_i:E\to V_{\lm_i}$ and the corresponding
automorphic trilinear functional $l^{aut}_{i}$ on $E\otimes
V_{\lm_i}^*$ defined by $l^{aut}_{i}(v\otimes u\otimes w)=\langle
r_\Dl(\nu_E(u\otimes v)),\bar \nu_i(w)\rangle$ (we identified
$\bar V_{\lm_i}$ with the smooth part of $V^*_{\lm_i}$). Such a
functional is clearly $G$-invariant, and hence we can invoke the
uniqueness principle for trilinear functionals from Section
\ref{triple}.

To this end, we fix a model trilinear functional
$l^{mod}_{\lm_i}=l^{mod}_{E\otimes V_{\lm_i}^*}$ (see Section
\ref{mod-functional} and the formula \eqref{K-circle} below; for a
detailed discussion, see \cite{BR3}) and the corresponding
intertwining model map $T_{\lm_i}=T_{\lm_i}^{mod}:E\to V_{\lm_i}$.
This gives rise to the coefficient of proportionality which we
denote by $a(\lm_i)=a_{\nu_E\otimes\nu_i}$ (somewhat abusing
notations by suppressing the dependence on $\nu_E$ and $\nu_i$)
such that $$T^{aut}_i=\nu_i^*(r_\Dl)=a(\lm_i)\cdot T_{\lm_i}\ .$$

The integral $p_\ck(f)=\int_\ck f(k)d\mu_\ck$ defines the period
map $p_\ck: C^\8(X)\to\bc$. Note that $p_\ck$ vanishes on
representations with no non-zero $K'$-invariant vectors, e.g., on
representations of discrete series.

We have the basic relation
$$d_{\Dl\ck}=(r_\Dl)_*(p_\ck)\ .$$
This means that for any $w\in E$, we have $d_{\Dl\ck}(w)=\int_\ck
r_\Dl(\nu_E(w))d\mu_\ck\ .$ We also have the following spectral
decomposition for any $w\in E$:
\begin{eqnarray}\label{K-diag-spec}
r_\Dl(w)=\sum_{\nu_i}\nu_i^*(r_\Dl(w))\ ,
\end{eqnarray} where $\nu_i$ runs through the fixed decomposition
of $L^2(X)$ into irreducible components.

We apply the functional $p_\ck$ to each term in
\eqref{K-diag-spec} and invoke the uniqueness principle for
$K'$-invariant functionals on irreducible representations
$V_{\lm_i}$ (i.e., that $d_{\lm_i}^{aut}=\beta(\lm_i)\cdot
d_{\lm_i}$; see Section \ref{P-n}). This, together with the
Fourier expansion \eqref{Fourier-d-K-w}, imply two different
expansions for the functional $d_{\Dl\ck}$: one which is {\it
``geometric''} (i.e., the Fourier expansion \eqref{Fourier-d-K-w}
along the orbit $\ck$) and another one which is {\it ``spectral''}
(i.e., induced by the trilinear invariant functionals and the
expansion \eqref{K-diag-spec}).

Namely, we have
\begin{align}\label{rs-K}
\sum_n|b_n(\nu)|^2 \hat
w(n,-n)=d_{\Dl\ck}(w)=\sum_{\lm_i}a(\lm_i)\beta(\lm_i)\cdot
d_{\lm_i}(T_{\lm_i}(w))\ ,
\end{align}
where $\hat w(n,-n)=\langle w,e'_n\otimes \bar e'_{-n}\rangle_E$
for any $w\in E$, with $\{e'_n\}$ a basis of $K'$-types in $V$ and
$\{\bar e'_n\}$ the conjugate basis in $\bar V$.

This is our substitute for the Rankin-Selberg formula in the
anisotropic case.

To make this formula explicit, we describe the model trilinear
functional in the circle model of representations $V=V_\tau$,
$\bar V=V_{-\tau}$ and $V_{\lm_i}$. We assume for simplicity that
$\tau\in i\br$ (i.e., $V$ is a representation of the principal
series) and that there is no exceptional spectrum for the lattice
$\G$ (i.e., that $\lm_i\in i\br$ for all $i>0$, and hence
$V^*_\lm\simeq V_{-\lm}$). The general case could be treated
analogously.

To simplify formulas, we make the following remark. The formula
\eqref{rs-K}  appeals only to automorphic representations of $G$
and a choice of a (non-trivial) connected compact subgroup
$K'\subset G$ (i.e., the choice of another compact subgroup $K$ we
made in Section \ref{K-geometry} is irrelevant). Since there is no
{\it preferred} compact subgroup in $G$, we may assume without
loss of generality that $K'=PSO(2)$ is the standard connected
compact subgroup of $G$.

It was shown in \cite{BR3} that in the circle model of class one
representations, the kernel of $l^{mod}_{E\otimes V_{-\lm}}$ on
the space $E\otimes V_{-\lm}\simeq C^\8_{even}(S^1\times S^1\times
S^1)$  is given by the following function in three variables
$\theta,\ \theta',\ \theta''\in S^1$
\begin{equation}\begin{split}
\label{K-circle}
K_{\tau,-\tau,\lm}(\theta,\theta',\theta'')=&\\
|\sin(\theta-\theta')|^{\frac{-1-\lm}{2}}&|\sin
(\theta-\theta'')|^{\frac{-1-2\tau+\lm}{2}}|\sin
(\theta'-\theta'')|^{\frac{-1+2\tau+\lm}{2}}\ .
\end{split}
\end{equation} This also defines the kernel of the map $T_\lm:E\to V_\lm$ via
the relation
\begin{align*} \langle
T_\lm(w),v\rangle_{V_\lm}=\frac{1}{(2\pi)^3}\int_{(S^1)^3}w(\theta,\theta')v(\theta')
K_{\tau,-\tau,\lm}(\theta,\theta',\theta'')d\theta
d\theta'd\theta''\ . \end{align*}

We have $d_{\lm}(T_\lm(w))=\langle
T_\lm(w),e_0\rangle_{V_\lm}=\frac{1}{(2\pi)^3}\int
w(\theta,\theta')K_{\tau,-\tau,\lm}(\theta,\theta',\theta'')d\theta
d\theta'd\theta''$ for any vector $w\in C^\8(S^1\times S^1)$. It
is clear from the formula \eqref{rs-K} that we can assume without
loss of generality that the vector $w\in E$ is $\Dl K$-invariant.
Such a vector $w$ can be described by a function of one variable;
namely, we set $w(\theta,\theta')= u(c)$ for $u \in C^\8(S^1)$ and
$c = (\theta - \theta')/2$. We have then $\hat w(n,-n)=\hat
u(n)=\frac{1}{2\pi}\int_{S^1} u(c)e^{-inc}dc$, i.e., $\hat u$ is
the Fourier transform of $u$.

We introduce a new kernel
\begin{eqnarray}\label{k-lm-def}
k_\lm(c)=k_{\tau,\lm} \left({\scriptstyle \frac{\theta -
\theta'}{2}}\right)=\frac{1}{2\pi}\int_{S^1}
K_{\tau,-\tau,\lm}(\theta,\theta',\theta'')d\theta''
\end{eqnarray} and the corresponding integral transform
\begin{eqnarray}\label{sperical-transf}
u^\sharp(\lm)=u^\sharp_\tau(\lm)=\frac{1}{(2\pi)^2}\int_{S^1}
u(c)k_{\lm}(c)dc\ ,
\end{eqnarray} suppressing the dependence on $\tau$ as we have fixed
the Maass form $\phi_\tau$. The transform is clearly defined for
any smooth function $u\in C^\8(S^1)$, at least for $\lm\in i\br$.
In fact, it could be defined for all $\lm\in \bc$, by means of
analytic continuation, but we will not need this.

Note that $k_\lm$ is the average of the kernel
$K_{\tau,-\tau,\lm}$ with respect to the action of $\Dl K$, or, in
other terms, is the pullback of the $K$-invariant vector $e_0\in
V_\lm$ under the map $T_\lm^*$, i.e., $k_\lm=T_\lm^*(e_0)\in E^*$.
We also note that the contribution in \eqref{rs-K} coming from the
trivial representation (i.e., $\lm=1$) is equal to
$u(0)=\frac{vol(\ck)}{vol(X)^\haf}\cdot  u(0)$ under our
normalization of measures $vol(X)=vol(\ck)=1$.

The formula \eqref{rs-K} then takes the form
\begin{eqnarray}\label{rs-K-u}
\sum_n|b_n(\nu)|^2 \hat u(n)=  u(0)+
\sum_{\lm_i\not=1}a(\lm_i)\beta(\lm_i)\cdot u^\sharp(\lm_i)\ .
\end{eqnarray}
This formula is an anisotropic analog of the Rankin-Selberg
formula \eqref{RS3} for the unipotent Fourier coefficients of
Maass forms. We finish the proof of Theorem \ref{thm-2}. \qed

\subsection{Remarks}\label{rem-sph}
1. The kernel function $k_\lm$ is not an elementary function,
unlike in the case of the unipotent Fourier coefficients where its
analog is given by $|x-y|^{-\haf-s}$. This is related to the fact
that the $N$-invariant distribution $\dl_s$ on $V_s$ is also
$\chi$-equivariant under the action of the full Borel subgroup
$B=AN$ for an appropriate character $\chi$ of $B$ which is trivial
on $N$. The space of $(B,\chi)$-equivariant distributions on $E$
is one-dimensional for a generic $\chi$. This is due to the fact
that $B$ has one open orbit for the diagonal action on the space
$\br\times \br$ and the vector space $E$ is modelled in the space
of smooth functions on this space. It is easy to write then a
non-zero $B$-equivariant functional on $E$ by an essentially
algebraic formula. We do not have a similar phenomenon for a
maximal compact subgroup of $G$. We will obtain however, an
elementary formula for leading terms in the asymptotic expansion
of $k_\lm$ as $|\lm|\to\8$ (see Appendix \ref{k-assymp}).

2. For Hecke-Maass forms, the proportionality coefficient $a(s)$
in the Rankin-Selberg formula \eqref{RS2} for the unipotent
Fourier coefficients gives the Rankin-Selberg $L$-function (after
multiplication by $\zeta(2s)$). In the anisotropic case we do not
know how to express the coefficient $a(\lm_i)$ in terms of an
appropriate $L$-function. It is known that the value of
$|a(\lm_i)|^2$ is related to the special value of the triple
$L$-function (see \cite{W}), but not the coefficient itself.  The
same is true for the coefficient $\beta(\lm_i)$ where in special
cases $|\beta(\lm_i)|^2$ is related to certain automorphic
$L$-function (see \cite{Wa}, \cite{JN}). There still might be a
way to normalize the product $a(\lm_i)\beta(\lm_i)$ in a canonical
way. We hope to return to this subject elsewhere.

3. For a non-uniform lattice $\G$ , the proof we gave above leads
to the following formula analogous to \eqref{rs-K-u} which
includes the contribution from the Eisenstein series. Namely, for
the similarly defined coefficients $a_k(s)$ and $\beta_k(s)$
corresponding to the Eisenstein series contribution, we have
\begin{eqnarray*}
 \sum_n|b_n(\nu)|^2 \hat u(n)=u(0)+
\sum_{\lm_i\not=1}a(\lm_i)\beta(\lm_i)\cdot u^\sharp(\lm_i)
+\sum_{\rm cusps}\int_{i\br^+} a_k(s)\beta_k(s)\cdot
u^\sharp(s)ds\ .\end{eqnarray*}

\subsection{Bounds for spherical Fourier coefficients.
Proof of Theorem \ref{thm-3}}\label{proof-thm-3} We follow the
strategy of Section \ref{proof-F-u}. We construct a $\Dl
K$-invariant test vector $w\in V\otimes \bar V$, i.e., a function
$u\in C^\8(S^1)$, such that when substituted into the
Rankin-Selberg formula \eqref{rs-K-u} will produce a weight $\hat
u$ which is not too small for a given $n$, $|n|\to\8$. We then
have to estimate the spectral density of such a vector, i.e., the
transform $u^\sharp$. In fact, as in Section \ref{proof-F-u} we
take a function which produces a weighted sum of the coefficients
$|b_k(\phi)|^2$ for $k$ in a short range depending on $n$ and such
that its transform $u^\sharp$ is spread over a relatively short
range of $\lm$'s. For such test vectors $w$ we give an essentially
sharp bound for the value of $d_{\Dl\mathcal{K}}(w)$.

We have the following technical

\begin{lem}{lem2}For any
integer $N$ and a real number $T$ such that  $N\geq T\geq 1$,
there exists a smooth function $u_{N,T}\in C^\8(S^1)$ satisfying
the following bounds:
\begin{enumerate}
    \item $|u_{N,T}(0)|\leq \al T$,
    \item $\hat u_{N,T}(k)\geq 0$ for all $k$,
    \item $\hat u_{N,T}(k)\geq 1$ for all $k$ satisfying $|k-N|\leq T$,
    \item $|u_{N,T}^\sharp(\lm)|\leq \al
    T|N|^{-\haf}(1+|\lm|)^{-\haf}+\al T(1+|\lm|^{-5/2})$ for $|\lm|\leq N/T$,
    $\lm\in i\br$,
    \item $|u_{N,T}^\sharp(\lm)|\leq \al
    T(1+|\lm|)^{-5/2}$ for $|\lm|\geq N/T$,
\end{enumerate} for some fixed constant $\al>0$ independent of $N$ and $T$.
\end{lem}

The proof of this Lemma is given in Appendix A. We construct the
corresponding function $u_{N,T}$ by considering a function
$u_{N,T}(c)=T^2\cdot e^{-iNc}\cdot\left(\psi_T*\psi_T'\right)(c)$,
where $*$ denotes the convolution in $C^\8(S^1)$ and
$\psi_T'(c)=\bar\psi_T(-c)$. Here $\psi_T(c)=\psi\left(T\cdot
c\right)$, where $\psi\in C^\8(S^1)$ is a {\it fixed} smooth
function supported in a small neighborhood of $1\in S^1$, and
$T\cdot c$ means the obvious scaling-up of the angle parameter in
$S^1$. Functions $u_{N,T}$ obviously satisfy bounds $(1)$-$(3)$
and the verification of $(4)$-$(5)$ is reduced to a routine
application of the stationary phase method (similar to our
computations in \cite{BR4}). These bounds are similar to bounds in
Lemma \ref{lem1} for the unipotent Fourier coefficients. There are
two differences though. First, in the corresponding bound in $(4)$
there is a factor of $(1+|\lm|)^{-\haf}$ in the first term. This
constitutes an important difference between $K$-invariant and
$N$-invariant functionals on the representation $E$. The
additional term $T(1+|\lm|^{-\frac{5}{2}})$ comes from the
estimate of the remainder in the stationary phase method and could
be improved further (although it would not make a difference in
what follows). The second (minor) difference is that the integral
transform $^\flat$ is elementary (i.e., the Mellin transform)
while the integral transform $^\sharp$ has its kernel given by a
non-elementary function (essentially by the hypergeometric
function). This slightly complicates computations.

\begin{rem}{remark-app}We would like to point out that it
is absolutely essential that in the integral
$u_{N,T}^\sharp(\lm)=\int_{S^1} u_{N,T}(c)k_\lm(c)\ dc$ the
support of the function $u_{N,T}$ does not contain points
$\pm\pi/4,\ \pm 3\pi/4$. Otherwise the phase in the above
oscillating integral possess {\it degenerate} critical points at
these values of $c$ for $N/T\asymp|\lm|$. The presence of
degenerate critical points change drastically the behavior of the
corresponding transform $^\sharp$. In particular, for
$n\asymp|\lm|$, the $\sharp$-transform of a {\it pure tensor}
$e_n\otimes\bar e_{-n}$ does not satisfy the bound $(4)$ in the
above lemma. Namely, $(e_n\otimes\bar e_{-n})^\sharp(\lm)$ have a
sharp peak for $n\asymp|\lm|$ of the order of
$|\lm|^{-\frac{5}{6}}$ (as oppose to $|\lm|\inv$). This phenomenon
is the starting point for the proof of the subconvexity bound for
the triple $L$-function given in \cite{BR4}. In the present paper,
we choose test vectors to vanish in a neighborhood of these
degenerate points $ \pm\pi/4,\ \pm 3\pi/4$ in the model
realization $ V\otimes \bar V\simeq C^\8(S^1\times S^1)$. This
allows us to avoid the more delicate analysis of degenerate
critical points. Note that our test vectors are not given by a
finite combination of {\it pure tensors} of $K$-types.
\end{rem}

We return to the proof of the theorem. In the proof we will use
two bounds for the coefficients $a(\lm_i)$ and $\beta(\lm_i)$.
Namely, it was shown in \cite{BR3} that
\begin{eqnarray}\label{convex-in-proof}
\sum\limits_{A\leq|\lm_i|\leq 2A}|a(\lm_i)|^2\leq a A^2\ ,
\end{eqnarray}  for
any $A\geq 1$ and some explicit $a>0$. The second  bound we will
need is the bound
\begin{eqnarray}\label{convex-in-proof-2}
\label{Horm}\sum_{A\leq|\lm_i|\leq 2A}|\beta(\lm_i)|^2\leq b A^2\
,\end{eqnarray} valid for any $A\geq 1$ and some $b$. In a
disguise this is the classical bound of L. H\"{o}rmander \cite{Ho}
for the average value at a point for eigenfunctions of the
Laplace-Beltrami operator on a compact Riemannian manifold (e.g.,
$\Dl$ on $Y$). This follows from the normalization
$|\beta(\lm_i)|^2=|\phi'_{\lm_i}(x_0)|^2$ we have chosen in
\eqref{b-0-nomalize} for $K'$-invariant eigenfunctions. In fact,
the bound \eqref{Horm} is standard in the theory of the Selberg
trace formula (see \cite{Iw}) and also can be easily deduced from
simple geometric considerations of \cite{BR3}. We note that both
bounds are consistent with the convexity bounds for relevant
$L$-functions.

We also assume for simplicity that there is no exceptional
spectrum (i.e., $\lm_i\in i\br$ for $i>0$). A general case could
be treated analogously by extending Lemma \ref{lem2} to cover
$\lm\in (-1,1)$ or as in \cite{BR4} where we treated exceptional
spectrum using simple considerations based on the Sobolev
restriction theorem.

We plug a test function satisfying bounds $(1)-(5)$, Lemma
\ref{lem2}, into the Rankin-Selberg formula \eqref{rs-K-u}. From
the Cauchy-Schwartz inequality, bounds \eqref{convex-in-proof} and
\eqref{Horm}, and summation by parts, we obtain

\begin{equation*}
\begin{split}
\sum_{|k-N|\leq T}|b_k(\nu)|^2\leq \sum_k|b_k(\nu)|^2\ \hat
  u_{N,T}(k)= u_{N,T}(0)
  +\sum_{\lm_i\not=1}a(\lm_i) \beta(\lm_i)&\ u_{N,T}^\sharp(\lm_i)\leq  \\
\leq \al T+\sum_{|\lm_i|\leq N/T}\al
T|N|^{-\haf}(1+|\lm_i|)^{-\haf}|a(\lm_i)\beta(\lm_i)|&\ +\\
+\sum_{\lm_i\not=1}\al
T(1+|\lm_i|)^{-\frac{5}{2}}|a(\lm_i)\beta(\lm_i)|& \leq\\
\leq \al T+\al T|N|^{-\haf}\cdot\sum_{|\lm_i|\leq
N/T}(1+|\lm_i|)^{-\haf}\left(|a(\lm_i)|^2+
  |\beta(\lm_i)|^2\right)&\ +\\
+\
T\sum_{\lm_i\not=1}(1+|\lm_i|)^{-\frac{5}{2}}\left(|a(\lm_i)|^2+
  |\beta(\lm_i)|^2\right)&\leq\\
  \leq \al T+CT|N|^{-\haf}\left(\frac{N}{T}\right)^{\frac{3}{2}+\eps}+DT=
  c'T+CT^{-\haf-\eps}|N|^{1+\eps}\ &,
\end{split}
\end{equation*} for any $\eps>0$ and some constants $c',\ C,\ D>0$.

Setting $T=N^{\frac{2}{3}}$, we obtain $ \sum\limits_{|k-N|\leq
N^{2/3}}|b_k(\nu)|^2\leq A_\eps N^{\frac{2}{3}+\eps}$ for any
$\eps>0$. \qed
\begin{remark}\label{rem-sph-lind}
One expects that bounds $|a(\lm_i)|\ll |\lm_i|^\eps$ and
$|\beta(\lm_i)|\ll |\lm_i|^\eps$ hold for any $\eps>0$. In special
cases this would be consistent with the Lindel\"{o}ff conjecture
for the corresponding $L$-functions. This however, will not have
the similar effect on the bound in Theorem \ref{thm-3} for
spherical Fourier coefficients $b_n(\phi_\tau)$. The reason for
such a discrepancy is that the spectral measure of the Eisenstein
series is much ``smaller" than that of the cuspidal spectrum.
Nevertheless, it is natural to conjecture that for general
$\G\subset\PGLR$ and a point $y_0\in Y$ the spherical Fourier
coefficients satisfy the bound $|b_n(\phi_\tau)|\ll|n|^\eps$. For
a CM-point $y_0$ and a Hecke-Maass form, this would correspond to
the Lindel\"{o}ff conjecture for the special value of the
corresponding $L$-function via the Waldspurger formula
\eqref{Wald}.
\end{remark}
%\newpage
\appendix
\section{Asymptotic expansions}\label{apen-A}
\subsection{Asymptotic expansion for the kernel
$k_\lm$}\label{k-assymp} We set $c=\frac{\theta - \theta'}{2}$ and
consider the integral \eqref{k-lm-def}, Section \ref{proof-thm-3}:
\begin{eqnarray*}
k_\lm(c)&=&k_{\tau,\lm} \left({\scriptstyle \frac{\theta -
\theta'}{2}}\right)=\frac{1}{2\pi}\int_{S^1}
K_{\tau,-\tau,\lm}(\theta,\theta',\theta'')\ d\theta''=\nonumber\\
&=&\frac{1}{2\pi}\cdot|\sin(2c)|^{-\haf-\hlm}\cdot\int_{S^1}
|\sin(\theta''-c)|^{-\haf-\tau+\hlm} |\sin(\theta''+c)|^{-\haf+\tau+\hlm}d\theta''\\
&=&|\sin(2c)|^{-\haf-\hlm}K_{\lm,\tau}(c) \ ,
\end{eqnarray*} where the kernel $K_{{-\tau,\tau,\lm}}$ is as in
\eqref{K-circle} and we denoted by $K_{\lm,\tau}$ the integral
\begin{eqnarray}\label{int-red2-B}K_{\lm,\tau}(c)=\frac{1}{2\pi}\cdot\int_{S^1}
|\sin(t-c)|^{-\haf-\tau+\hlm} |\sin(t+c)|^{-\haf+\tau+\hlm}dt\ .
\end{eqnarray}

The kernel $K_{\lm,\tau}$ is not given by an elementary function.
We obtain an asymptotic formula for $K_{\lm,\tau}$ by applying the
stationary phase method to the integral \eqref{int-red2-B}. The
asymptotic formula we obtain is valid for a {\it fixed} $\tau$ and
is uniform for $\lm\in i\br$ and $c\not=0,\ \pm\pi/2,\ \pi$. We
denote the set of exceptional points by $S=\{0,\ \pm\pi/2,\
\pi\}\subset S^1$.  These points are singular because the
integrand in \eqref{int-red2-B} degenerates for these values of
$c$. We have the following

\begin{claim}{claim}There are constants $A$, $B$ and $C$ such that
for all $\lm\in i\br$ and $c\not\in S$,
\begin{eqnarray}\label{m+r}
K_{\lm,\tau}(c)=m_\lm(c)+m_\lm(c+\pi/2)+r_{\tau}(\lm,c)\ ,
\end{eqnarray}
where the main term $m_\lm(c)$ is a smooth function of $\lm$ and
$c$ ($c\not\in S$), and for $|\lm|\geq1$ is given by
\begin{eqnarray}\label{int-K-lm}
m_\lm(c)=|\lm|^{-\haf}\left(A+B|\lm|\inv+C|\lm|\inv\cos^2(c)\right)
\cdot|\sin(c)|^{\lm}\ .
\end{eqnarray} The remainder $r_{\tau}(\lm,c)$ satisfies the estimate
\begin{equation}\label{rem}
|r_{\tau}(\lm,c)|=
O\left((1+|\lm|)^{-5/2}+\left[1+\left|\ln(|\sin(c)\cos(c)|)\right|\right]
\cdot (1+|\lm|)^{-10}\right)\
\end{equation} with the implied constant in the $O$-term depending
on $\tau$ only.
\end{claim}

Hence we have
\begin{equation}
\begin{split}\label{k-lm-asym}k_\lm(c)=
|\sin(2c)|^{-\haf-\hlm}K_{\lm,\tau}(c)=&\\
M_\lm(c)+M_{\lm}(c+\pi/2)+
&|\sin(2c)|^{-\haf-\hlm}r_{\tau}(\lm,c)\ , \
\end{split}\end{equation} with
$M_\lm(c)=|\lm|^{-\haf}\left[A+B|\lm|\inv+C|\lm|\inv\cos^2(c)\right]\cdot
|\sin(2c)|^{-\haf}|\sin(c)|^{\hlm}|\cos(c)|^{-\hlm}$.

\subsection{Proof} The asymptotic expansion in the claim follows from the
stationary phase method applied to the integral
\eqref{int-red2-B}. We consider the asymptotic expansion
consisting of two terms and a reminder. Since all functions are
$\pi$-periodic we consider only the interval $c\in [0,\pi]$. For
$c\not\in S$, the phase of the oscillating kernel in the integral
\eqref{int-red2-B} has two non-degenerate critical points at $t=0$
and $t=\pi/2$. Hence, the asymptotic expansion is given by the sum
of two terms. It turns out that these terms have the form
$m_\lm(c)$ and $m_\lm(c+\pi/2)$ for the same function $m_\lm$.
Singularities of the amplitude at $t=c,\ \pi-c$ are responsible
for the logarithmic term in the remainder. For $|\lm|\to\8$, this
contribution from the singularities of the amplitude is of order
of $O((1+|\lm|)^{-k})$ for any $k>0$ due to the fast oscillation
of the phase at the same points. In fact, in \cite{BR4} we gave a
self-contained  treatment of such (and more complicated) integrals
based on the reduction to standard integrals and the use of the
van der Corput lemma. Here we show how one can deduce necessary
bounds from the stationary phase method.

Our computations are based on the following well-known form of the
two-term asymptotic in the stationary phase method (see \cite{Bo},
\cite{F}, \cite{St}). We also use the estimation of the
corresponding remainder.

Let $\phi$ and $f$ be smooth real valued functions on $S^1$. To
state the stationary phase formula, we assume that $\phi$ has a
unique non-degenerate critical point $t_0\in S^1$. We consider the
integral $I(\lm)=\int_{S^1}f(t)e^{\lm\phi(t)}dt$ for $\lm\in
i\br$. For $|\lm|\geq 1$, we have the following expansion
\begin{eqnarray}\label{stat-phase-formula}
I(\lm)=|\lm|^{-\haf}(C_0+C_1|\lm|\inv)e^{\lm\phi(t_0)}+r(\lm)\
,\end{eqnarray} where $C_0=(2\pi)^\haf e^{i\cdot{\rm
sign(\phi''(t_0))}\pi/4}|\phi''(t_0)|^{-\haf}f(t_0)$ and
\begin{equation*}
\begin{split}
 C_1=(\pi/2)^\haf\ e^{3i\cdot{\rm
sign(\phi''(t_0))}\pi/4}\ |\phi''(t_0)|^{-\frac{3}{2}}
 \cdot &\\ \cdot[f''-\phi^{(3)}f'/\phi'' - &
\phi^{(4)}f/4\phi''+5(\phi^{(3)})^2f/12(\phi'')^2]_{t=t_0}\ ,
  \end{split}
\end{equation*}
and the remainder satisfies $r(\lm)=O((1+|\lm|)^{-5/2})$. The
constant in the $O$-term is bounded for $\phi$ and $f$ in a
bounded, with respect to natural semi-norms, set in $C^\8(S^1)$.
If $\phi$ has a number of isolated non-degenerate critical points
then the asymptotic is given by the sum over these points of the
corresponding contributions.

For $|\lm|<1$, we have the trivial bound: $|I(\lm)|\leq \int|f|\
dt$.

\subsubsection{Leading terms}
We apply these formulas to compute leading terms
in the asymptotic expansion of the integral \eqref{int-red2-B}. We
set
$$\phi(t)=\ln|\sin(t-c)|+\ln|\sin(t+c)|\ \ {\rm and}
\ \ f(t)=|\sin(t-c)|^{-\haf-\tau} |\sin(t+c)|^{-\haf+\tau}\ .$$ We
have $\phi'(t)=\sin(2t)/\sin(t-c)\sin(t+c)$ and hence the phase
$\phi$ has two critical points $t=0$ and $t=\pi/2$, assuming that
$c\not=0,\ \pi/2,\ \pi$.

A straightforward computation gives for $t=0$,
$$\phi''(0)=-2\sin^{-2}(c),\ \phi^{(3)}(0)=0,\
\phi^{(4)}(0)=-4(1+2\cos^2(c))/\sin^4(c)$$ and
$f(0)=|\sin(c)|\inv,\ f''(0)=|\sin(c)|^{-3}(1+4\tau^2\cos^2(c))\
,$ and similarly for $t=\pi/2$,
$$\phi''(\pi/2)=-2\cos^{-2}(c),\ \phi^{(3)}(\pi/2)=0,\
\phi^{(4)}(\pi/2)=-4(1+2\sin^2(c))/\cos^4(c)$$ and
$f(\pi/2)=|\cos(c)|\inv,\
f''(\pi/2)=|\cos(c)|^{-3}(1+4\tau^2\sin^2(c))\ .$

Plugging this into \eqref{stat-phase-formula} we see that for
$c\not=0,\ \pm\pi/2,\ \pi$,
\begin{eqnarray}\label{asymp-K-lm}
K_{\lm,\tau}(c)=m_\lm(c)+m_\lm(c+\pi/2)+ r(\lm,c)\ ,
\end{eqnarray} where $m_\lm(c)=|\lm|^{-\haf}\left(A+B|\lm|\inv+C|\lm|\inv\cos^2(c)\right)
\cdot|\sin(c)|^{\lm}$ with some explicit constants $A,\ B,\ C$.

\subsubsection{The remainder}\label{rem-app} We need to estimate the remainder
$r(\lm,c)=r_{\tau}(\lm,c)$ as $c$ keeps away from the singular set
$S\subset S^1$. We claim that
\begin{align*}|r_{\tau}(\lm,c)|=O\left((1+|\lm|)^{-5/2}+
\left[1+|\ln(|\sin(c)\cos(c)|)|\right]\cdot
(1+|\lm|)^{-10}\right)\ . \end{align*} To see this we first note
that the contribution to the remainder coming from the integration
over any {\it fixed} interval which does not include singularities
of the phase and of the amplitude (i.e., points $t=\pm c,\ \pi\pm
c$) is of order $O((1+|\lm|)^{-5/2})$ (with the constant in the
$O$-term depending on the proximity of the interval to these
singular points and on $\tau$).

The analysis of the integral is identical for all points in $S$,
hence we treat only the case of $0<c\leq\frac{\pi}{4}$.

We use the appropriate partition of unity in order to separate
different behavior of the kernel. Let $i_0, i_c, i_{-c},
i_{\pi/2}, i_{-\pi/2}, i_{\pi-c}, i_{\pi+c}, i_\pi\subset
S^1=\br/2\pi\bz$ be the following collection of closed {\it
overlapping} intervals: $i_0=[-3c/4,3c/4],\ i_c=[c/4,\pi/3],\
i_{-c}=[-\pi/3,-c/4],\ i_\pi=\pi+i_0,\ i_{\pi-c}=\pi+i_{-c},\
i_{\pi+c}=\pi+i_c$ and similarly  $i_{\pi/2}=[\pi/4,3\pi/4],\
i_{-\pi/2}=[-\pi/4, -3\pi/4]$. Let
$1=\chi_0+\chi_c+\chi_{-c}+\chi_{\pi/2}+\chi_{-\pi/2}+\chi_{\pi-c}+\chi_{\pi+c}+
\chi_\pi$ be the corresponding partition of the unity on $S^1$
separating singular points $\pm c,\ \pi\pm c$, from the stationary
points $0,\ \pi,\ \pm\pi/2$ (i.e., $supp(\chi_0)\subset i_0$
etc.). For $\chi_i$ as above, we denote by $I^i_{\lm,\tau}$ the
corresponding  integral
\begin{eqnarray*} I^{i}_{\lm,\tau}(c)=\frac{1}{2\pi}\cdot\int_{S^1}
|\sin(t-c)|^{-\haf-\tau+\hlm}
|\sin(t+c)|^{-\haf+\tau+\hlm}\chi_{i}(t)dt\ .
\end{eqnarray*}
We have $K_{\lm,\tau}(c)=\sum I^{i}_{\lm,\tau}(c)$. Due to the
symmetry, it is enough to deal with the integrals $I^0_{\lm,\tau},
I^c_{\lm,\tau}, I^{\pi/2}_{\lm,\tau}$.

The integral $I^{\pi/2}_{\lm,\tau}(c)$ falls under the standard
stationary phase method and hence the remainder in the two-term
asymptotic is of order of $O([1+|\lm|]^{-5/2})$ (with the constant
in the $O$-term depending on $\tau$). In fact the behavior of the
integral $I^{0}_{\lm,\tau}(c)$ is similar. This could be seen
easily by scaling-up the variable $t$ by $c$. In particular,
integrals $I^{\pm\pi/2}_{\lm,\tau}(c)$ and
$I^{0,\pi}_{\lm,\tau}(c)$ give rise to the leading terms in the
asymptotic in Claim \ref{claim} and the remainder which is of the
order of $O([1+|\lm|]^{-5/2})$.

The behavior of the remaining integral  $I^c_{\lm,\tau}$ is
similar to the well-known Beta type integral of the form
$B(\lm,\chi)=\int
|x+1|^{-\haf-\tau+\hlm}|x-1|^{-\haf+\tau+\hlm}\chi(x)\ dx$ for a
compactly supported smooth function $\chi$, vanishing in a {\it
neighborhood} of the stationary point (i.e., near $x=0$). (In
fact, in \cite{BR4} we showed how to reduce the integral
$I^c_{\lm,\tau}$ to the integral $B(\lm,\chi)$ using an
appropriate change of variable.) It follows from integration by
parts that for such $\chi$, the integral $B(\lm,\chi)$ is of the
order of $O([1+|\lm|]^{-k})$ for any $k>0$. The similar analysis
is applicable to the integral $I^c_{\lm,\tau}$.

For $|\lm|\leq 1$, the integral $I^c_{\lm,\tau}$ is trivially of
the order of $O(|\ln(|\sin(c)\cos(c)|)|)$.

For $|\lm|>1$ and small $c$,  we consider two integrals
\begin{eqnarray*}
J^{-}_{\lm,\tau}(c)=\frac{1}{2\pi}\cdot\int_{c/4}^{c}
|\sin(t-c)|^{-\haf-\tau+\hlm}
|\sin(t+c)|^{-\haf+\tau+\hlm}\chi_{c}(t)dt\ ,\\
J^{+}_{\lm,\tau}(c)=\frac{1}{2\pi}\cdot\int_{c}^{\pi/3}
|\sin(t-c)|^{-\haf-\tau+\hlm}
|\sin(t+c)|^{-\haf+\tau+\hlm}\chi_{c}(t)dt\ .
\end{eqnarray*}
Scaling-up the variable $t$ by $c$, we see that the integral
$J^{-}_{\lm,\tau}(c)$ transforms into an integral of the form
$\int_{\frac{1}{4}}^1|f(x)|^{-\haf+\lm/2}\psi(x)$, where
$f(x)=f_c(x)$ is a monotone smooth function with derivatives
$f^{(n)}$ bounded for $x\in[1/4, 1]$, uniformly in $c$, and
satisfying $f(1)=0$ and $f'(x)\geq 1$ for $x\in[1/4,1]$.  The
function $\psi$ is a smooth function also with, uniformly in $c$,
bounded derivatives (depending on $\tau$). Hence integration by
parts implies that the integral $J^{-}_{\lm,\tau}(c)$ is of the
order of $O([1+|\lm|]^{-k})$ for any $k>0$ with the constant in
the $O$-term independent of $c$.

Scaling-up the variable $t$ by $c$, we see that the integral
$J^{+}_{\lm,\tau}(c)$ transforms into an integral of the form
$\int_1^{\frac{\pi}{3c}}|(x-1)\cdot g(x)|^{-\haf+\lm/2}\psi(x)$,
where $g$ is a monotone smooth function with, uniformly in $c$,
bounded derivatives, satisfying $g(1)=1$ and $g'(x)\geq 1$ for
$x\in[1, \pi/3c]$, and $\psi$ is a smooth function with, uniformly
in $c$, bounded derivatives (depending on $\tau$). To estimate
such an integral, one breaks the interval $[1,\pi/3c]$ into
$|\ln(c)|$ dyadic intervals. On each such interval we have the
bound as above of the order of $O([1+|\lm|]^{-k})$ for any $k>0$
with the constant in the $O$-term independent of $c$. Hence the
integral $J^{+}_{\lm,\tau}(c)$ is of the order of
$O(|\ln(c)|\cdot[1+|\lm|]^{-k})$ for any $k>0$. \qed

\subsection{Proof of Lemma \ref{lem2}}\label{prf-lm2}
We have to analyze the integral transform given by
$u_{N,T}^\sharp(\lm)=\int u_{N,T}(c)k_\lm(c)dc$, where
$u_{N,T}(c)=T^2\cdot e^{-iNc}\cdot\left(\psi_T*\psi'_T\right)(c)$
with the parameters $N>T\geq 1$, $N\in\bz$, $T\in \br$. Here
$\psi\in C^\8(S^1)$ is a fixed smooth function of a support in a
small interval containing $1\in S^1$, and $\psi_T(c)=\psi(T\cdot
c)$, $\psi'_T(c)=\bar\psi_T(-c)$.

To analyze the asymptotic of $u_{N,T}^\sharp(\lm)$,  we consider
the following model integral
\begin{eqnarray*}
I(\lm,N,T)=T\int
e^{-iNc}|\sin(2c)|^{-\haf}|\sin(c)|^{\hlm}|\cos(c)|^{-\hlm}\chi(Tc)\
dc\ ,
\end{eqnarray*} where $\chi$ is a  smooth function with $supp(\chi)\subset[-1,1]$.

On the basis of the asymptotic expansion \eqref{k-lm-asym} for the
kernel $k_\lm$, we see that $u_{N,T}^\sharp(\lm)$ is of the order
of $I(\lm,N,T)\cdot (1+|\lm|)^{-\haf}+ O(T(1+|\lm|)^{-5/2})$.

We claim that for $|\lm|\leq N/T$, we have
$|I(\lm,N,T)|=O(TN^{-\haf})$ and for $|\lm|>N/T$, we have
$|I(\lm,N,T)|=O(|\lm|^{-k})$ for any $k>0$. These bounds imply the
claim in Lemma \ref{lem2}.

To obtain desired bounds for $I(\lm,N,T)$, we appeal again to the
stationary phase method.

 Namely, scaling-up by $T$ the variable
$c$ in the integral $I(\lm,N,T)$, we arrive at the integral
\begin{eqnarray*}
I_1(\lm,N,T)=\int e^{-i\frac{N}{T}
t}|\sin({\scriptstyle\frac{2}{T}}t)|^{-\haf}
|\tan({\scriptstyle\frac{t}{T}})|^{\hlm}\chi(t)dt\ .
\end{eqnarray*}
It is easy to see that for $|\lm|\leq 1$, this integral is of the
same order as the integral
$T^\haf\int|t|^{-\haf+\hlm}e^{-i\frac{N}{T} t}\chi(t)dt $, which
is of the order of $O(TN^{-\haf})$. For $1<|\lm|\leq N/T$, the
phase function in the integral $I_1$ has unique non-degenerate
critical point and the contribution from the singularities of the
amplitude is negligible. Hence, arguing as in Section
\ref{rem-app}, we see that the integral  $I_1$  is of the order of
$O(T N^{-\haf})$.

In fact, both cases follow immediately from the van der Corput
lemma (see \cite{BR4} for similar bounds).

For $|\lm|>N/T$, the phase function has no critical points and
hence we have  $|I_1|=O([1+|\lm|]^{-k})$ for any $k>0$.\qed

%%%%%%%%%%%%%%%%%%%%%%%%%%%%%%%%%%%%%%%%%%%%%%%


\begin{thebibliography}{BR2}
\bibitem[Be]{Be} J. Bernstein, Eisenstein series, lecture notes.
Park City, Utah (2004).

\bibitem[BR1]{BR1} J. Bernstein\ and\ A. Reznikov,
Analytic continuation of representations and estimates of
automorphic forms, Ann. of Math. (2) {\bf 150} (1999), no.~1,
329--352. MR1715328 (2001h:11053), arXiv: math.RT/9907202.

\bibitem[BR2]{BR2} \bysame, Sobolev norms of automorphic functionals,
Int. Math. Res. Not. {\bf 2002}, no.~40, 2155--2174. MR1930758
(2003h:11058)

\bibitem[BR3]{BR3}\bysame, Estimates of automorphic functions,
Mosc. Math. J. {\bf 4} (2004), no.~1, 19--37, 310. MR2074982 (2005f:11097),
arXiv: math.RT/0305351.

\bibitem[BR4]{BR4} \bysame, Subconvexity of
triple $L$-functions, preprint,  arXiv: math.NT/0608555 (2006).

\bibitem[B]{B} A. Borel, {\it Automorphic forms on ${\rm SL}\sb 2({\bf R})$},
Cambridge Univ. Press, Cambridge, 1997. MR1482800 (98j:11028)

\bibitem[Bo]{Bo} V. A. Borovikov, {\it Uniform stationary phase method},
IEE, London, 1994. MR1325462 (96b:58111)

\bibitem[Bu]{Bu} D. Bump, The Rankin-Selberg method: an introduction and survey,
in {\it Automorphic representations,
$L$-functions and applications: progress and prospects},
41--73, de Gruyter, Berlin. MR2192819 (2006k:11097)

\bibitem[BM]{BM} R. Bruggeman, Y. Motohashi, A new approach to the spectral theory
of the fourth moment of the Riemann zeta-function. J. Reine Angew.
Math. 579 (2005), 75--114. MR2124019 (2005i:11108)

\bibitem[F]{F} M. Fedoruk, Asymptotic methods in analysis. {\it Analysis. I},
Encyclopaedia Math. Sci., 13, Springer, Berlin, 1989. MR1042759 (90j:00027)


\bibitem[G5]{G5} I. M. Gelfand, M. I. Graev\ and\ N. Ya. Vilenkin,
{\it Generalized functions. Vol. 5},
Academic Press, New York, 1966. MR0207913 (34 \#7726)

\bibitem[G6]{GGPS} I. M. Gelfand, M. I. Graev\ and\
I. I. Pyatetskii-Shapiro, {\it Representation theory and
automorphic functions}, Saunders, Philadelphia, 1969. Reprint of
the 1969 edition. Generalized Functions, 6. Academic Press, Inc.,
Boston, MA, 1990. MR0233772 (38 \#2093).


\bibitem[Go]{Go1} A. Good, Cusp forms and eigenfunctions of the Laplacian,
Math. Ann. {\bf 255} (1981), no.~4, 523--548. MR0618183
%(82i:10029)

\bibitem[Gr]{Gr} B. H. Gross, Some applications of Gelfand pairs
to number theory, Bull. Amer. Math. Soc. (N.S.) {\bf 24} (1991),
no.~2, 277--301. MR1074028 (91i:11055)

\bibitem[He]{He} S. Helgason, {\it Groups and geometric analysis},
Amer. Math. Soc., Providence, RI, 2000. MR1790156 (2001h:22001)

\bibitem[Ho]{Ho} L. H\"ormander, {\it The analysis of linear partial
differential operators. IV}, Springer, Berlin, 1985. MR0781537 (87d:35002b)

\bibitem[Iw]{Iw} H. Iwaniec, {\it Spectral methods of automorphic forms},
Amer. Math. Soc., Providence, RI, 2002. MR1942691 (2003k:11085)

\bibitem[IS]{IS} H. Iwaniec\ and\ P. Sarnak,
Perspectives on the analytic theory of $L$-functions, Geom. Funct. Anal.
{\bf 2000}, Special Volume, Part II, 705--741. MR1826269 (2002b:11117)

\bibitem[JN]{JN} H. Jacquet\ and\ N. Chen, Positivity of quadratic base change $L$-functions,
Bull. Soc. Math. France {\bf 129} (2001), no.~1, 33--90. MR1871978 (2003b:11048)

\bibitem[KSa]{KSa} H. Kim, P. Sarnak,  Appendix to: H. H. Kim,
Functoriality for the exterior square of ${\rm GL}\sb 4$ and the
symmetric fourth of ${\rm GL}\sb 2$, J. Amer. Math. Soc. {\bf 16}
(2003), no.~1, 139--183. MR1937203 (2003k:11083)

\bibitem[KS]{KS} B. Kr\"otz\ and\ R. J. Stanton,
Holomorphic extensions of representations. I. Automorphic
functions, Ann. of Math. (2) {\bf 159} (2004), no.~2, 641--724.
MR2081437 (2005f:22018)

\bibitem[Ku]{Ku} T. Kubota, {\it Elementary theory of Eisenstein series},
Kodansha, Tokyo, 1973. MR0429749 (55 \#2759)

\bibitem[Kz]{Kz} N. V. Kuznetsov, Sums of Kloosterman sums and the eighth
power moment of the Riemann zeta-function, in {\it Number theory
and related topics (Bombay, 1988)}, 57--117, Tata Inst. Fund.
Res., Bombay. MR1441327 (98c:11085)

\bibitem[L]{L} J. B. Lewis, Eigenfunctions on symmetric spaces with
distribution-valued boundary forms, J. Funct. Anal. {\bf 29}
(1978), no.~3, 287--307. MR0512246 (80f:43020)

\bibitem[M]{M} H. Maass, \"Uber eine neue Art von nichtanalytischen automorphen
Funktionen und die Bestimmung Dirichletscher Reihen durch
Funktionalgleichungen, Math. Ann. {\bf 121} (1949), 141--183. MR0031519 (11,163c)

\bibitem[MW]{MW} K. Martin,  D. Whitehouse, Central $L$-values and toric periods for $GL(2)$,
preprint (2007).

\bibitem[Mo1]{Mo1} Y. Motohashi, {\it Spectral theory of the Riemann
zeta-function}, Cambridge Univ. Press, Cambridge, 1997. MR1489236 (99f:11109)

\bibitem[Mo2]{Mo2} \bysame, A note on the meanvalue of the zeta
and $L$-functions, preprint, arXiv: math.NT/0401085.

\bibitem[O]{O} A. I. Oksak, Trilinear Lorentz invariant forms,
Comm. Math. Phys. {\bf 29} (1973), 189--217. MR0340478 (49 \#5231)

\bibitem[PS]{PS} Y. N. Petridis\ and\ P. Sarnak, Quantum unique
ergodicity for ${\rm SL}\sb 2(\mathcal{O})\backslash\bold H\sp 3$
and estimates for $L$-functions, J. Evol. Equ. {\bf 1} (2001),
no.~3, 277--290. MR1861223 (2003a:11060)

\bibitem[Pr]{Pr} D. Prasad, Trilinear forms for representations
of ${\rm GL}(2)$ and local $\epsilon$-factors, Compositio Math.
{\bf 75} (1990), no.~1, 1--46. MR1059954 (91i:22023)

\bibitem[Ra]{Ra}  R. A. Rankin, Contributions to the theory of Ramanujan's
function $\tau(n)$ and similar arithmetical functions.  II. The
order of the Fourier coefficients of integral modular forms, Proc.
Cambridge Philos. Soc. {\bf 35} (1939), 357--372. MR0000411
(1,69d)

\bibitem[R]{R} A. Reznikov, Norms of geodesic restrictions
for eigenfunctions on hyperbolic surfaces and representation
theory, preprint (2004). arXiv: math.AP/0403437.

\bibitem[Sa]{Sa3} P. Sarnak, Integrals of products of eigenfunctions,
Internat. Math. Res. Notices {\bf 1994}, no.~6, 251--261.
MR1277052
%(95i:11039)

\bibitem[Se]{Se} A. Selberg, On the estimation of Fourier coefficients
of modular forms, in {\it Proc. Sympos. Pure Math., Vol. VIII},
1--15, Amer. Math. Soc., Providence, R.I. MR0182610 (32 \#93)

\bibitem[St]{St} E. M. Stein, {\it Harmonic analysis:
real-variable methods, orthogonality, and oscillatory integrals},
Princeton Univ. Press, Princeton, NJ, 1993. MR1232192 (95c:42002)


\bibitem[V]{V} A. Venkatesh, Sparse equidistribution problems, period bounds,
and subconvexity, preprint (2005). arXiv: math.NT/0506224.

\bibitem[Wa]{Wa} J.-L. Waldspurger, Sur les valeurs de certaines fonctions
$L$ automorphes en leur centre de sym\'etrie,
Compositio Math. {\bf 54} (1985), no.~2, 173--242. MR0783511 (87g:11061b)

\bibitem[W]{W} T. Watson, Thesis, Princeton, 2001.


\end{thebibliography}
\end{document}